\documentclass[final]{siamltex}
\usepackage[T1]{fontenc}
\usepackage[latin9]{inputenc}
\usepackage{geometry}
\usepackage{xcolor}
\usepackage{pdfcolmk}
\usepackage{amsbsy}
\usepackage{amstext}
\usepackage{amssymb}
\usepackage{graphicx}
\usepackage{esint}
\PassOptionsToPackage{normalem}{ulem}
\usepackage{ulem}
\usepackage[unicode=true,pdfusetitle,
 bookmarks=true,bookmarksnumbered=false,bookmarksopen=false,
 breaklinks=false,pdfborder={0 0 1},backref=section,colorlinks=false]
 {hyperref}

\makeatletter

\providecolor{lyxadded}{rgb}{1,0,0}
\providecolor{lyxdeleted}{rgb}{0,0,1}

\newcommand\eqref[1]{(\ref{#1})}

\@ifundefined{date}{}{\date{}}
\usepackage{ae,aecompl}

\usepackage{color}

\usepackage{amsfonts}
\usepackage{latexsym}\usepackage{color}

\usepackage{bm}

\newcommand{\sfrac}[2]{\mathchoice
  {\kern0em\raise.5ex\hbox{\the\scriptfont0 #1}\kern-.15em/
   \kern-.15em\lower.25ex\hbox{\the\scriptfont0 #2}}
  {\kern0em\raise.5ex\hbox{\the\scriptfont0 #1}\kern-.15em/
   \kern-.15em\lower.25ex\hbox{\the\scriptfont0 #2}}
  {\kern0em\raise.5ex\hbox{\the\scriptscriptfont0 #1}\kern-.2em/
   \kern-.15em\lower.25ex\hbox{\the\scriptscriptfont0 #2}}
  {#1\!/#2}}

\makeatother

\begin{document}

\title{Efficient Variable-Coefficient Finite-Volume Stokes Solvers}

\author{Mingchao Cai \thanks{Courant Institute of Mathematical Sciences, New York University, New York, NY 10012, E-mail: cmchao2005@gmail.com} \and Andy Nonaka \thanks{Center for Computational Sciences and Engineering, Lawrence Berkeley National Laboratory, Berkeley, CA 94720, E-mail: ajnonaka@lbl.gov} \and John B. Bell \thanks{Center for Computational Sciences and Engineering, E-mail: jbbell@lbl.gov} \and Boyce E. Griffith \thanks{Leon H. Charney Division of Cardiology, Department of Medicine, New York University School of Medicine, NY, and Courant Institute of Mathematical Sciences, E-mail: griffith@cims.nyu.edu} \and Aleksandar Donev \thanks{Corresponding author, Courant Institute of Mathematical Sciences, E-mail: donev@courant.nyu.edu}}
\maketitle
\begin{abstract}
We investigate several robust preconditioners for solving the saddle-point
linear systems that arise from spatial discretization of unsteady
and steady variable-coefficient Stokes equations on a uniform staggered
grid. Building on the success of using the classical projection method
as a preconditioner for the coupled velocity-pressure system {[}\emph{B.
E. Griffith, J. Comp. Phys., 228 (2009), pp. 7565\textendash{}7595}{]},
as well as established techniques for steady and unsteady Stokes flow
in the finite-element literature, we construct preconditioners that
employ independent generalized Helmholtz and Poisson solvers for the
velocity and pressure subproblems. We demonstrate that only a single
cycle of a standard geometric multigrid algorithm serves as an effective
inexact solver for each of these subproblems. Contrary to traditional
wisdom, we find that the Stokes problem can be solved nearly as efficiently
as the independent pressure and velocity subproblems, making the overall
cost of solving the Stokes system comparable to the cost of classical
projection or fractional step methods for incompressible flow, even
for steady flow and in the presence of large density and viscosity
contrasts. Two of the five preconditioners considered here are found
to be robust to GMRES restarts and to increasing problem size, making
them suitable for large-scale problems. Our work opens many possibilities
for constructing novel unsplit temporal integrators for finite-volume
spatial discretizations of the equations of low Mach and incompressible
flow dynamics.
\end{abstract}

\global\long\def\V#1{\boldsymbol{#1}}
 \global\long\def\M#1{\boldsymbol{#1}}
 \global\long\def\Set#1{\boldsymbol{#1}}

\global\long\def\D#1{\Delta#1}
 \global\long\def\d#1{\delta#1}

\global\long\def\norm#1{\left\Vert #1\right\Vert }
 \global\long\def\abs#1{\left|#1\right|}

\global\long\def\grad{\M{\nabla}}
\global\long\def\avv#1{\langle#1\rangle}
\global\long\def\av#1{\left\langle #1\right\rangle }

\global\long\def\Fb{\boldsymbol{F}}
 \global\long\def\gb{\boldsymbol{g}}
 \global\long\def\mb{\boldsymbol{m}}
 \global\long\def\bA{\boldsymbol{\mathcal{A}}}
 \global\long\def\bB{\boldsymbol{\mathcal{B}}}
 \global\long\def\bI{\boldsymbol{\mathcal{I}}}
 \global\long\def\bR{\boldsymbol{\mathcal{R}}}
\global\long\def\bzero{\boldsymbol{0}}

\global\long\def\Ab{\boldsymbol{A}}
 \global\long\def\Db{\boldsymbol{D}}
 \global\long\def\fb{\boldsymbol{f}}
 \global\long\def\Gb{\boldsymbol{G}}
 \global\long\def\Ib{\boldsymbol{I}}
 \global\long\def\Lb{\boldsymbol{L}}
 \global\long\def\Mb{\boldsymbol{M}}
 \global\long\def\Pb{\boldsymbol{P}}
 \global\long\def\ub{\boldsymbol{u}}

\global\long\def\deltab{\boldsymbol{\delta}}
 \global\long\def\phib{\boldsymbol{\phi}}
 \global\long\def\taub{\boldsymbol{\tau}}
 \global\long\def\psib{\boldsymbol{\psi}}

\global\long\def\half{\frac{1}{2}}
 \global\long\def\myhalf{\sfrac{1}{2}}
 \global\long\def\threehalf{\sfrac{3}{2}}

\global\long\def\IN{\boldsymbol{N}}
 \global\long\def\IZ{\boldsymbol{Z}}
 \global\long\def\IR{\boldsymbol{R}}
 \global\long\def\IC{\boldsymbol{C}}

\section{Introduction}

Many numerical methods for solving the time-dependent (unsteady) incompressible
\cite{bellColellaGlaz:1989,almgren-iamr,NonProjection_Griffith,LargeViscosityJump}
or low Mach number \cite{LowMachAdaptive,LowMachExplicit} equations
require the solution of a linear unsteady Stokes flow subproblem.
The linear steady Stokes problem is of particular interest for low
Reynolds number flows \cite{OlshanskiiReusken,VariableViscosity_May}
or flow in viscous boundary layers. In this work, we investigate efficient
linear solvers for the unsteady and steady Stokes equations in the
presence of variable density and viscosity. Specifically, we consider
the coupled velocity-pressure Stokes system \cite{Turekbook,Elman_book}
\begin{equation}
\left\{ \begin{array}{ll}
{\rho}\ub_{t}+\grad p=\grad\cdot{\taub}({\ub})+\fb,\\
\grad\cdot\ub=g,
\end{array}\right.\label{Stokes_model}
\end{equation}
where ${\rho\left(\V r\right)}$ is the density, ${\ub\left(\V r,t\right)}$
is the velocity, $p\left(\V r,t\right)$ is the pressure, $\V f\left(\V r,t\right)$
is a force density, and $\taub(\ub)$ is the viscous stress tensor.
A nonzero velocity-divergence $g\left(\V r,t\right)$ arises, for
example, in low Mach number models because of compositional or temperature
variations \cite{LowMachAdaptive}. The viscous stress $\taub(\ub)$
is ${\mu}\grad\ub$ for constant viscosity incompressible flow, ${\mu}\left[\grad\ub+(\grad\ub)^{T}\right]$
when $g=0$ (incompressible flow), and ${\mu}\left[\grad\ub+(\grad\ub)^{T}\right]+({\gamma}-\frac{2}{3}{\mu})(\grad\cdot\ub)\Ib$
when $g\ne0$, where ${\mu\left(\V r,t\right)}$ is the shear viscosity
and ${\gamma\left(\V r,t\right)}$ is the bulk viscosity. When the
inertial term is neglected, $\rho\ub_{t}=0$, (\ref{Stokes_model})
reduces to the time-independent (steady) Stokes equations. In this
work we consider periodic boundary conditions and physical boundary
conditions that involve velocity only, notably no-slip and free-slip
physical boundaries%
\footnote{When the normal component of velocity is specified on the whole boundary
of the computational domain $\Omega$, a compatibility condition $\int_{\partial\Omega}\V u\cdot\V n\, dS=\int_{\Omega}g\, d\V r$
needs to be imposed.%
}.

Spatial discretization of (\ref{Stokes_model}) can be carried out
using standard finite-volume or finite-element techniques. Applying
the backward Euler scheme to solve the spatially-discretized equations
with time step size $\D t$, gives the following discrete system for
the velocity $\ub^{n+1}$ and the pressure $p^{n+1}$ at the end of
time step $n$, 
\begin{equation}
\left\{ \begin{array}{ll}
{\rho}\left(\frac{\ub^{n+1}-\ub^{n}}{\Delta t}\right)+\grad p^{n+1}=\grad\cdot{\taub}\left({\ub}^{n+1}\right)+\fb^{n+1},\\
\grad\cdot\ub^{n+1}=g^{n+1},
\end{array}\right.\label{Stokes_eqn}
\end{equation}
where $\fb^{n+1}$ contains external forcing terms such as gravity
and any explicitly-handled terms such as, for example, advection.
Similar linear systems are obtained with other implicit and semi-implicit
temporal discretizations \cite{bellColellaGlaz:1989,almgren-iamr,NonProjection_Griffith}.
In the limit $\rho/\D t\rightarrow0$, the system (\ref{Stokes_eqn})
reduces to the steady Stokes equations. Here we will assume that the
spatial discretization is stable, more precisely, that the Stokes
system (\ref{Stokes_eqn}) is ``uniformly solvable'' as the spatial
discretization becomes finer, i.e., that a suitable measure of the
condition number of the Schur complement of (\ref{Stokes_eqn}) remains
bounded as the grid spacing $h\rightarrow0$. In the context of finite-element
methods, this is typically implied by the well-known $\inf-\sup$
or Ladyženskaja-Babuška-Brezzi (LBB) condition. Here we employ the
classical staggered-grid \cite{HarWel65} discretization on a uniform
grid, which can be thought of as a rectangular analog of the lowest-order
Raviart-Thomas element and is known to be a stable discretization
\cite{StaggeredGrid_Stability,Olshanskii}. We expect that it will
be relatively straightforward to generalize the preconditioners developed
here to recently-developed adaptive mesh staggered schemes \cite{IBAMR_HeartValve,VariableViscosity_May}.
Note, however, that collocated finite-volume discretizations of the
Navier-Stokes equations do not provide a stable discretization, which
motivates the development of approximate-projection methods \cite{almgrenBellSzymczak:1996}.

Historically, there have been significant differences in the treatment
of (\ref{Stokes_eqn}) in the finite-volume and finite-element literature.
In the finite-element literature, there is a long history of numerical
methods for solving the Stokes equations, especially in the time-independent
(steady) context \cite{Turekbook,Elman_book}. By contrast, in the
context of high-resolution finite-volume methods, the dominant paradigm
has been to use a splitting (fractional-step) or projection method
\cite{Chorin68,ProjectionMethods_Minion} to separate the pressure
and velocity updates. In part, this choice has been motivated by the
target applications, which are often high Reynolds number, or even
inviscid, flows. In the inviscid limit, the splitting error associated
with projection methods vanishes, and for sufficiently large Reynolds
number flows the time step size dictated by advective stability constraints
makes the splitting error relatively small. At the same time, the
preference for splitting methods stems, in large part, from the perception
that solving the saddle-point problem (\ref{Stokes_eqn}) is much
more difficult than solving the pressure and velocity subproblems;
to quote the authors of Ref. \cite{ProjectionMethods_Minion}, ``Spatially
discretized versions of the coupled Eqs. ... are cumbersome to solve
directly.'' In fact, one of the first second-order projection methods
\cite{bellColellaGlaz:1989} was developed by starting with a Crank-Nicolson
variant of (\ref{Stokes_eqn}) and then approximating the resulting
Stokes system using a velocity-pressure splitting that was motivated
by the perceived difficulty in solving the coupled system.

Fractional-step approaches, however, suffer from several significant
shortcomings. It is well-known, for example, that the splitting introduces
a commutator error that leads to the appearance of spurious or ``parasitic''
modes \cite{GaugeIncompressible_E,ProjectionMethods_Minion} in the
presence of physical boundaries. Furthermore, it is generally not
possible to impose the true boundary conditions of the Stokes system
in a fractional-step scheme; instead, \textquotedbl{}artificial\textquotedbl{}
boundary conditions must be imposed in the velocity and pressure subsystems.
This has motivated the construction of methods that \emph{approximately}
solve the Stokes system (\ref{Stokes_eqn}) using block-triangular
factorizations \cite{BlockTriangularProjection} similar to those
employed here to construct preconditioners for iterative methods that
solve (\ref{Stokes_eqn}) \emph{exactly}. This is crucial at small
Reynolds numbers because the splitting error becomes larger as viscous
effects become more dominant, and projection or approximate factorization
methods do not apply in the steady Stokes regime for problems with
physical boundary conditions.

Recognizing these problems, one of us investigated the use of projection-like
methods as preconditioners for a Krylov method for solving the coupled
system (\ref{Stokes_eqn}) \cite{NonProjection_Griffith}. It was
found that, contrary to traditional finite-volume wisdom and in agreement
with extensive experience in the finite-element context, the saddle-point
problem (\ref{Stokes_eqn}) can be efficiently solved using standard
multigrid techniques for the velocity and pressure subproblems, for
a broad range of parameters. Here we improve and generalize the preconditioners
developed in Ref. \cite{NonProjection_Griffith} to account for variable
density and variable viscosity, as well as to robustly handle small
or zero Reynolds number (steady) flows. Our primary motivation for
this work is the development of semi-implicit integrators for the
low Mach number equations of fluctuating hydrodynamics for multicomponent
fluid mixtures \cite{LowMachExplicit}. For these applications it
is important to treat the viscosity implicitly (including the limit
of steady Stokes flow \cite{DiffusionJSTAT}) due to the large separation
of time scales between momentum and mass diffusion (i.e., large Schmidt
number). In order to properly include thermal fluctuations with implicit
viscous handling it is also necessary to use a coupled Stokes formulation
instead of split (projection method) approaches \cite{LLNS_Staggered,DFDB}.

The preconditioners we investigate here numerically are drawn from
the large finite-element literature on Stokes solvers \cite{Benzi_JCP,Benzi,Elman,Elman_book,Elman2007,Ipsen,Kay,MardalWinther2004,MardalWinther2011,Murphy1,OlshanskiiReusken,OlshanskiiMG,Turekbook}.
We will not attempt to review the extensive finite-element literature
on preconditioners for Stokes flow here; instead, we will point out
the similarities and differences with prior work for each of the preconditioners
that we study. When necessary, we generalize the existing preconditioners
to finite Reynolds numbers (unsteady Stokes equations) and to variable
density and variable viscosity problems. We investigate several alternative
preconditioners that solve the velocity and pressure subproblems in
different orders. In the finite-element context, variable-viscosity
steady Stokes solvers based on several of the preconditioners we investigate
here have been developed by several groups \cite{VariableViscosity_FEM,VariableViscosity_FEM2},
and have already been successfully scaled to massively-parallel architectures
and very difficult large-contrast geophysical problems. In the finite-volume
context, the work most closely related to our work is Ref. \cite{LargeViscosityJump},
which focuses on steady Stokes flow in the presence of large viscosity
contrast (i.e., discontinuities) for geodynamic applications. Notably,
both our work and the work presented in Ref. \cite{LargeViscosityJump}
are based on a staggered finite-volume discretization and geometric
multigrid solvers.

Our primary contribution in this work is that we investigate in detail
the computational performance of a collection of four standard preconditioners
over a broad range of parameters in the context of a specific but
very efficient (in terms of number of degrees of freedom per grid
cell) finite-volume discretization. By carefully designing and optimizing
the parameters for all of the key components of the solvers, ranging
from the geometric multigrid smoothers to the restart frequency of
the GMRES solver, we construct a complete solver that can readily
be employed to construct novel unsplit temporal integrators for finite-volume
(conservative) spatial discretizations of the equations of low Mach
and incompressible flow dynamics. Importantly, these unsplit schemes
can use the same building blocks (e.g., geometric multigrid solvers
and high-resolution advection techniques) and achieve a similar computational
complexity as traditional projection methods, as we will demonstrate
in future work.

The preconditioners that we investigate are built using two crucial
subsolvers. The first of these is a linear solver for the inviscid
problem
\begin{equation}
\left\{ \begin{array}{ll}
{\rho}\left(\frac{\ub^{n+1}-\ub^{n}}{\Delta t}\right)+\grad p^{n+1}=\fb^{n+1},\\
\grad\cdot\ub^{n+1}=g^{n+1},
\end{array}\right.\label{eq:pressure_subproblem}
\end{equation}
the solution of which requires solving a density-weighted pressure
Poisson equation
\[
-\D t\,\grad\cdot\left(\rho^{-1}\grad p^{n+1}\right)=g^{n+1}-\grad\cdot\left(\ub^{n}+\rho^{-1}\fb^{n+1}\D t\right).
\]
For the staggered-grid finite-volume discretization we employ here,
this Poisson problem can efficiently be solved using standard geometric
multigrid techniques \cite{almgren-iamr}. The second subsolver required
by the preconditioners is a linear solver for the unconstrained variable-coefficient
velocity equation,
\begin{equation}
{\rho}\left(\frac{\ub^{n+1}-\ub^{n}}{\Delta t}\right)=\grad\cdot{\taub}\left({\ub}^{n+1}\right)+\fb^{n+1}.\label{eq:velocity_subproblem}
\end{equation}
Note that both (\ref{eq:pressure_subproblem}) and (\ref{eq:velocity_subproblem})
use the same boundary conditions for velocity as the coupled problem,
and that natural boundary conditions are required for the pressure
when solving (\ref{eq:pressure_subproblem}) on a staggered grid.
For constant viscosity incompressible flow $\grad\cdot{\taub}\left({\ub}\right)=\mu\grad^{2}{\ub}$
and therefore (\ref{eq:velocity_subproblem}) is a system of $d$
uncoupled Helmholtz equations, where $d$ is the dimensionality. These
can be solved efficiently using standard geometric multigrid techniques.
For variable viscosity flows or when $g\neq0$ the different components
of velocity are coupled. Here we develop an effective geometric multigrid
method for solving (\ref{eq:velocity_subproblem}) that generalizes
the classical red-black coloring smoother for the scalar Poisson equation.
Since the solution of either (\ref{eq:pressure_subproblem}) or (\ref{eq:velocity_subproblem})
is itself a costly iterative process, it is crucial that the preconditioners
require only approximate subsolvers. More precisely, preconditioning
should only require the application of linear operators that are spectrally-equivalent
\cite{Elman_book} to the exact solution operators for (\ref{eq:pressure_subproblem})
or (\ref{eq:velocity_subproblem}). Here we use one or a few cycles
of geometric multigrid as approximate solvers for these subproblems.
The preconditioners investigated in this work can be easily generalized
to other spatial discretizations and boundary condition types by simply
modifying the approximate subsolvers for (\ref{eq:pressure_subproblem})
and (\ref{eq:velocity_subproblem}). For example, boundary conditions
that couple pressure and viscous stress can be handled by imposing
approximate boundary conditions for the subsolvers. In Ref. \cite{NonProjection_Griffith},
at physical boundaries on which normal tractions (normal components
of the stress tensor) are prescribed, Neumann conditions are imposed
on the normal velocity component when solving (\ref{eq:velocity_subproblem})
and Dirichlet conditions are imposed for the pressure when solving
(\ref{eq:pressure_subproblem}). For adaptively-refined meshes \cite{IBAMR_HeartValve,VariableViscosity_May},
multilevel geometric multigrid techniques can be used to solve the
pressure and velocity subproblems \cite{almgren-iamr,IBAMR_HeartValve}.

The organization of this paper is as follows. In section \ref{sec:Preconditioners},
we introduce several preconditioners based on approximating the inverse
of the Schur complement. In Section \ref{sec:Staggered} we specialize
to a particular staggered-grid second-order finite-volume discretization
and give details of our numerical implementation. In Section \ref{sec:Results}
we perform a detailed study of the efficiency and robustness of the
various preconditioners, and select the optimal values for several
algorithmic parameters. Finally, we offer some conclusions in Section
\ref{sec:Conclusions}, and then give several technical derivations
in an extensive Appendix.

\section{\label{sec:Preconditioners}Preconditioners}

In this section we present several preconditioners for solving the
saddle-point linear system (\ref{Stokes_eqn}) that arises after spatio-temporal
discretization of (\ref{Stokes_model}). Much of the discussion presented
here has already appeared scattered through many diverse works in
the literature; for the benefit of the reader we provide a condensed
but complete summary of the key derivations. For increased generality,
we write this system in the form, 
\begin{equation}
\V M\left(\begin{array}{c}
\V x_{\ub}\\
\V x_{p}
\end{array}\right)=\left(\begin{array}{cc}
\Ab & \Gb\\
-\Db & \V 0
\end{array}\right)\left(\begin{array}{c}
\V x_{\ub}\\
\V x_{p}
\end{array}\right)=\left(\begin{array}{c}
\V b_{\ub}\\
\V b_{p}
\end{array}\right),\label{Saddle_sys}
\end{equation}
where $(\V x_{\ub},\,\V x_{p})^{T}$ denote the velocity and pressure
degrees of freedom, $(\V b_{\ub},\,\V b_{p})^{T}$ are the velocity
and pressure right hand sides, $\Db$ denotes a discrete divergence
operator, and $\V G$ is a discrete gradient operator. Note that for
the staggered-grid discretization that we describe in Section \ref{sec:Staggered},
the gradient and divergence operators are negative adjoints of each
other for periodic, no-slip, and free-slip boundary conditions, $\V G=(-\Db)^{*}$,
where star denotes adjoint, making $\V M=\V M^{\star}$ a self-adjoint
matrix. Here the linear velocity operator $\Ab=\theta{\V{\rho}}-\V L_{\V{\mu}}$
combines inertial and viscous effects, where $\theta$ is a parameter
that is zero for steady Stokes flow, and $\theta\sim\D t^{-1}$ for
unsteady flow. The operator ${\V{\rho}}$ is a mass density matrix
(distinct from the standard finite element mass matrix), such that
$\V{\rho}\V x_{\ub}$ is a spatially-discrete (conserved) momentum
field. The viscous operator is denoted with $\V L_{\V{\mu}}$, with
$\V L_{\V{\mu}}\ub$ being a spatial discretization of $\grad\cdot\taub(\ub)$.

The saddle-point problem (\ref{Saddle_sys}) can formally be solved
by using the inverse of the Schur complement, 
\[
{\displaystyle \V S^{-1}=\left(-\V D\Ab^{-1}\Gb\right)^{-1},}
\]
to obtain the exact solution for the pressure, 
\begin{equation}
\V x_{\V p}=-\V S^{-1}({\Db}\Ab^{-1}\V b_{\ub}+\V b_{p}),\label{sol_exact_p}
\end{equation}
and for the velocity degrees of freedom,
\begin{equation}
\V x_{\ub}=\Ab^{-1}(\V b_{\ub}-\Gb\V x_{p})=\Ab^{-1}\V b_{\ub}+\Ab^{-1}\Gb\V S^{-1}(\Db\Ab^{-1}\V b_{\ub}+\V b_{p}).\label{sol_exact_u}
\end{equation}
These formal solutions are not useful in practice because the Schur
complement cannot be formed explicitly for large three-dimensional
grids, nor inverted efficiently. In Ref. \cite{LargeViscosityJump},
the authors investigate evaluating the action of $\V S^{-1}$ in (\ref{sol_exact_p})
by an outer Krylov solver, which itself relies on evaluating the action
of $\Ab^{-1}$ in an inner (nested) Krylov solver. We do not investigate
this approach here and instead focus on what the authors of Ref. \cite{LargeViscosityJump}
call the ``fully coupled preconditioned approach'', in which an
approximation of the Schur complement solution is used to construct
an effective preconditioner for a Krylov solver applied to the saddle-point
problem (\ref{Saddle_sys}). The key part in designing preconditioners
for (\ref{Saddle_sys}) is approximating the (inverse of the) Schur
complement, specifically, constructing an operator $\mathcal{S}^{-1}$
that is spectrally-equivalent to $\V S^{-1}$ \cite{Elman}.

To motivate the approximation of $\V S^{-1}$, let us consider the
case of constant viscosity $\mu_{0}$ and constant density $\rho_{0}$.
In this case $\Ab=\theta\rho_{0}\Ib-\mu_{0}\Lb$, where $\M I$ denotes
an identity matrix and $\V L$ is a discrete vector Laplacian operator,
constructed taking into account the imposed velocity boundary conditions.
We then have
\begin{equation}
\V S^{-1}=\left[-\Db\left(\theta\rho_{0}\V I-\mu_{0}\V L\right)^{-1}\Gb\right]^{-1}\approx\left[\left(-\V D\V G\right)\left(\theta\rho_{0}\V I-\mu_{0}\M L_{p}\right)^{-1}\right]^{-1}=-\theta\rho_{0}\M L_{p}^{-1}+\mu_{0}\V I,\label{Schur_commutApp}
\end{equation}
where $\M L_{p}=\Db\Gb$ denotes a scalar (pressure) discrete Laplacian
operator, and have assumed the commuting property $\M L\M G\approx\Gb\M L_{p}$,
which is an exact identity for the staggered grid discretization applied
to periodic systems. This approximation to the Schur complement inverse
has been used in the finite-element context in Ref. \cite{KayGresho}
and in the finite-volume approach in Ref. \cite{NonProjection_Griffith};
an in-depth discussion of the use of approximate commutators for constructing
preconditioners can be found in Ref. \cite{ApproximateCommutators}.

Here we generalize (\ref{Schur_commutApp}) to variable density and
viscosity through a simple construction. The basic idea is that the
first part of the Schur complement approximation, $\theta\rho_{0}\M L_{p}^{-1}$,
corresponds to the inviscid limit. For variable density, this term
becomes $\theta\V L_{\V{\rho}}^{-1}$, where 
\[
{\Lb}_{\V{\rho}}=\Db{\V{\rho}}^{-1}\Gb
\]
is a discretization of the density-weighted Poisson operator $\grad\cdot\rho^{-1}\grad$
that also appears in traditional variable-density projection methods
\cite{almgren-iamr}. Therefore, for variable-density, constant-viscosity
flow, $\grad\cdot\taub={\mu}_{0}\grad^{2}\ub$, and we employ the
approximation
\begin{equation}
\V S^{-1}\approx\mathcal{S}^{-1}=-\theta\V L_{\V{\rho}}^{-1}+\mu_{0}\V I.\label{eq:Schur_app_constVisc}
\end{equation}

The term $\mu_{0}\V I$ in (\ref{Schur_commutApp}) is an analogue
of the viscous operator $\V L_{\V{\mu}}$ that acts on pressure-like
degrees of freedom instead of velocity-like degrees of freedom. This
has to be constructed on a case-by-case basis, and in the constant
viscosity setting it corresponds to the viscous pressure-correction
term proposed by Brown, Cortez and Minion \cite{ProjectionMethods_Minion}
in the context of second-order projection methods. For incompressible
flow, $\taub(\ub)={\mu}\left[\grad\ub+(\grad\ub)^{T}\right]$, the
Fourier-space calculation described in Appendix \ref{AppendixFourier}
suggests replacing the term $\mu_{0}\V I$ with $2{\V{\mu}}$, where
${\V{\mu}}$ is a diagonal matrix of viscosities corresponding to
each pressure degree of freedom. This gives the Schur complement inverse
approximation
\begin{equation}
\V S^{-1}\approx\mathcal{S}^{-1}=-\theta{\Lb}_{\V{\rho}}^{-1}+2{\V{\mu}},\label{Schur_app_varDen_varVis}
\end{equation}
which is called the ``local viscosity'' preconditioner in Ref. \cite{LargeViscosityJump}.
Note however that the prefactor of two suggested by the analysis in
Appendix \ref{AppendixFourier} is not included in Eq. (36) in Ref.
\cite{LargeViscosityJump}. When bulk viscosity is included, $\taub(\ub)={\mu}\left[\grad\ub+(\grad\ub)^{T}\right]+({\gamma}-\frac{2}{3}{\mu}(\grad\cdot\ub))\Ib$,
we take 
\begin{equation}
{\displaystyle \V S^{-1}\approx\mathcal{S}^{-1}=-\theta\V L_{\V{\rho}}^{-1}+\left({\V{\gamma}}+\frac{4}{3}{\V{\mu}}\right),}\label{eq:Schur_app_bulk}
\end{equation}
where ${\V{\gamma}}$ is the diagonal matrix of bulk viscosities.
As we demonstrate in Appendix \ref{AppendixFourier}, these approximations
are exact for periodic systems if the density and viscosity are constant.
In all other cases they are approximations that are expected to be
good in regions far from boundaries where the coefficients do not
vary significantly. Our numerical experiments support this intuition.

We have investigated the alternative approximations 
\begin{equation}
\V S^{-1}\approx-\theta{\Lb}_{\V{\rho}}^{-1}-\V L_{\V{\rho}}^{-1}\Db{\V{\rho}}^{-1}\V L_{\V{\mu}}{\V{\rho}}^{-1}\Gb\V L_{\V{\rho}}^{-1},\label{eq:BFBt}
\end{equation}
as well as
\[
\V S^{-1}\approx-\theta{\Lb}_{\V{\rho}}^{-1}-\V L_{p}^{-1}\left(\Db\V L_{\V{\mu}}\Gb\right)\V L_{p}^{-1},
\]
which is similar to the so-called BFBt preconditioner of Elman \cite{Elman}
in the steady-state case, and which is also investigated in Ref. \cite{LargeViscosityJump}.
These approximations utilize the velocity boundary conditions since
they involve the viscous operator $\V L_{\V{\mu}}$, unlike the pressure-space
viscous operator in (\ref{Schur_app_varDen_varVis}) which does not
make use of the velocity boundary conditions. We have observed similar
behavior for the more expensive approximation (\ref{eq:BFBt}) as
with the simpler and significantly more efficient approximation (\ref{Schur_app_varDen_varVis}).
We therefore do not investigate BFBt-type preconditioners in this
work.

As explained in Appendix \ref{sec:AnalysisEigenvalues}, the spectrum
of the preconditioned operators for the preconditioners we consider
next is determined by the spectrum of $\mathcal{S}^{-1}\V S$. In
that appendix, we demonstrate with a combination of analytical techniques
and numerical computation that this operator has a very clustered
spectrum even in the presence of non-trivial boundary conditions and
large variations in viscosity.

\subsection{Projection Preconditioner}

In the first preconditioner we consider, which we will denote with
$\V P_{1}$, we use one step of the classical projection method \cite{Chorin68,bellColellaGlaz:1989,NonProjection_Griffith}
as a preconditioner. In $\V P_{1}$, we use (\ref{sol_exact_p}) to
estimate the pressure, and make a commuting assumption in (\ref{sol_exact_u}),
\[
\Ab^{-1}\Gb\V S^{-1}=\Ab^{-1}\Gb(-\Db\Ab^{-1}\Gb)^{-1}\approx-\Ab^{-1}\Ab{\V{\rho}}^{-1}\Gb\Lb_{\V{\rho}}^{-1}=-{\V{\rho}}^{-1}\Gb\Lb_{\V{\rho}}^{-1},
\]
which gives the velocity estimate 
\begin{equation}
\V x_{\ub}\approx\Ab^{-1}\V b_{\ub}-{\V{\rho}}^{-1}\Gb\Lb_{\V{\rho}}^{-1}(\Db\Ab^{-1}\V b_{\ub}+\V b_{p}).\label{eq:velocity_projection}
\end{equation}
Note that this velocity estimate (\ref{eq:velocity_projection}) satisfies
the divergence condition exactly, $\Db\V x_{\ub}=-\V b_{p}$. More
precisely, $\V x_{\ub}$ is the $L_{2}$ projection of the unconstrained
velocity estimate $\Ab^{-1}\V b_{\ub}$ onto the divergence constraint.

In practical implementation, the exact subproblem solvers need to
be replaced by approximations. Specifically, ${\Ab}^{-1}$ is approximated
by the inexact velocity solver $\widetilde{\Ab}^{-1}$, $\V L_{\V{\rho}}^{-1}$
is implemented by the approximate pressure Poisson solver $\widetilde{\V L}_{\V{\rho}}^{-1}$,
and $\V S^{-1}$ is replaced by $\widetilde{\mathcal{S}}^{-1}$, which
is an approximation to the approximate Schur complement inverse $\mathcal{S}^{-1}$
given by (\ref{Schur_app_varDen_varVis}) for incompressible flow.
In summary, for the variable-coefficient Stokes problem, the projection
preconditioner $\V P_{1}$ is defined by the block factorization 
\begin{equation}
\V P_{1}^{-1}=\left(\begin{array}{cc}
\V I & {\V{\rho}}^{-1}\V G\widetilde{\V L}_{\V{\rho}}^{-1}\\
\V 0 & \widetilde{\mathcal{S}}^{-1}
\end{array}\right)\left(\begin{array}{cc}
\V I & \V 0\\
-\V D & -\V I
\end{array}\right)\left(\begin{array}{cc}
\widetilde{\Ab}^{-1} & \V 0\\
\V 0 & \V I
\end{array}\right).\label{P1_practical}
\end{equation}
This factorization clearly shows the main steps in the application
of the preconditioner. First, a velocity subproblem is solved inexactly
(right-most block) to compute $\V x_{\ub}^{*}=\widetilde{\Ab}^{-1}\V b_{\ub}$.
Second, $\V b_{c}={\Db}\V x_{\ub}^{*}+\V b_{p}$ is computed (middle
block). Third, a Poisson problem is solved approximately to compute
$\widetilde{\V L}_{\V{\rho}}^{-1}\V b_{c}$ and, lastly, the pressure
and velocity estimates are evaluated (first block). For constant-coefficient
periodic problems with exact subsolvers, the projection preconditioner
is an exact solver for the coupled Stokes equations since both (\ref{Schur_commutApp})
and (\ref{eq:velocity_projection}) are exact.

For the constant viscosity and density Stokes problem, a projection
preconditioner very similar to $\V P_{1}$ was first proposed by one
of us in Ref. \cite{NonProjection_Griffith}. In this work we generalize
the projection preconditioner to the case of variable viscosity and
density. Even in the constant-coefficient case, there is a small but
important difference between $\V P_{1}$ and the previous projection
preconditioner in Ref. \cite{NonProjection_Griffith}, which uses
the following approximation of the Schur complement inverse, 
\[
\V S^{-1}\approx\widetilde{\mathcal{S}}^{-1}=-\left(\theta\rho_{0}\V I-\mu_{0}\M L_{p}\right)\widetilde{\V L}_{p}^{-1},
\]
rather than the approximation (\ref{eq:Schur_app_constVisc}) used
here, $\widetilde{\mathcal{S}}^{-1}=-\theta\rho_{0}\widetilde{\V L}_{p}^{-1}+\mu_{0}\V I$,
which we have found to give a slightly more efficient solver. The
two approximations are identical when exact Poisson solvers are used,
$\widetilde{\V L}_{p}^{-1}=\M L_{p}^{-1}$, but not when an approximate
solver is employed.

\subsection{Lower Triangular Preconditioner}

For our second preconditioner, which we denote with $\V P_{2}$, we
use (\ref{sol_exact_p}) for the pressure estimate, but the velocity
estimate takes the simpler form 
\begin{equation}
\V x_{\ub}\approx{\Ab}^{-1}\V b_{\ub},\label{sol_u_inexact}
\end{equation}
which is obtained by discarding the second part in (\ref{sol_exact_u}).
If we further approximate the matrix inverses with inexact solves,
namely, replacing ${\Ab}^{-1}$ by $\widetilde{\Ab}^{-1}$, $\V L_{\V{\rho}}^{-1}$
by $\widetilde{\V L}_{\V{\rho}}^{-1}$, and $\V S^{-1}$ by $\widetilde{\mathcal{S}}^{-1}$,
the second preconditioner is given by the block factorization 
\begin{equation}
\V P_{2}^{-1}=\left(\begin{array}{cc}
\Ib & \V 0\\
\V 0 & -\widetilde{\mathcal{S}}^{-1}
\end{array}\right)\left(\begin{array}{cc}
\Ib & \V 0\\
\Db & \Ib
\end{array}\right)\left(\begin{array}{cc}
\widetilde{\Ab}^{-1} & \V 0\\
\V 0 & \Ib
\end{array}\right).\label{CC_Pre_lower}
\end{equation}

By combing all the terms in the right hand side of (\ref{CC_Pre_lower}),
we see that $\V P_{2}$ is actually an approximation of the inverse
of the lower triangular preconditioner previously studied by several
other groups \cite{CahouetChabard,Ipsen,Kay,MardalWinther2004,MardalWinther2011},
\begin{equation}
\V P_{2}^{-1}\approx\left(\begin{array}{cc}
\Ab & \V 0\\
-\V D & -\V S
\end{array}\right)^{-1}.\label{eq:P2_S}
\end{equation}
Notice that for steady Stokes flow, $\theta=0$, the application of
$\V P_{2}^{-1}$ does not require any pressure Poisson solvers, unlike
the projection preconditioner. Therefore, a single application of
$\V P_{2}^{-1}$ can be significantly less expensive computationally
than an application of $\V P_{1}^{-1}$. For unsteady flows $\V P_{1}$
and $\V P_{2}$ involve nearly the same operations and applying them
has similar computational cost.

\subsection{Upper Triangular Preconditioner}

Alternatively, one can assume $\Db\Ab^{-1}\V b_{\ub}\approx0$ to
obtain $\V x_{p}=-\V S^{-1}\V b_{p}$ and 
\[
{\displaystyle \V x_{\ub}={\Ab}^{-1}(\V b_{\ub}+\Gb\V S^{-1}\V b_{p}).}
\]
Replacing the exact solvers with inexact solvers, we obtain our third
preconditioner in block factorization form, 
\begin{equation}
\V P_{3}^{-1}=\left(\begin{array}{cc}
\widetilde{\Ab}^{-1} & \V 0\\
\V 0 & \Ib
\end{array}\right)\left(\begin{array}{cc}
\Ib & -\Gb\\
\V 0 & \Ib
\end{array}\right)\left(\begin{array}{cc}
\Ib & \V 0\\
\V 0 & -\widetilde{\mathcal{S}}^{-1}
\end{array}\right),\label{CC_Pre_upper}
\end{equation}
which is exactly the same as the ``fully coupled'' approach with
the ``local viscosity'' preconditioner studied in Ref. \cite{LargeViscosityJump}
and also the block-triangular preconditioner of Ref. \cite{VariableViscosity_FEM2},
generalized here to time-dependent problems. If we combine all the
terms in the right hand side of (\ref{CC_Pre_upper}), then we see
that $\V P_{3}$ is actually an approximation of the inverse of the
upper triangular preconditioner \cite{CahouetChabard,Ipsen,Kay,MardalWinther2004,MardalWinther2011},
\begin{equation}
\V P_{3}^{-1}\approx\left(\begin{array}{cc}
\Ab & \Gb\\
\V 0 & -\V S
\end{array}\right)^{-1}.\label{eq:P3_S}
\end{equation}
The computational cost of applying $\V P_{3}^{-1}$ is very similar
to that of applying $\V P_{2}^{-1}$.

\subsection{Other preconditioners}

In addition to the three main preconditioners (projection, lower and
upper triangular) we study here, we have investigated some other preconditioners.
The simplest Schur-complement based preconditioner one can construct
is the block diagonal preconditioner \cite{CahouetChabard,Kay,MardalWinther2004,MardalWinther2011}
\begin{equation}
\V P_{4}^{-1}=\left(\begin{array}{cc}
\widetilde{\Ab}^{-1} & \V 0\\
\V 0 & -\widetilde{\mathcal{S}}^{-1}
\end{array}\right).\label{CC_Pre_diagonal}
\end{equation}
This preconditioner has the lowest computational cost of all the preconditioners
per Krylov iteration, but it also yields the poorest approximation
to the exact solution (\ref{sol_exact_p},\ref{sol_exact_u}). Note,
however, that the use of a diagonal preconditioner can make the preconditioned
operator symmetric and thus allow for the use of more efficient (short-recurrence)
Krylov solvers such as MINRES. This is exploited in Ref. \cite{VariableViscosity_FEM}
to construct a robust and highly-scalable finite-element discretization
of the variable-viscosity steady Stokes equations, using a single
cycle of algebraic multigrid for a Laplacian approximation to $\Ab$
as an approximate velocity solver.

In Appendix \ref{sec:AnalysisEigenvalues}, we show that $\V P_{1}$,
$\V P_{2}$ and $\V P_{3}$ all give the same spectrum for the preconditioned
linear operator. It is also well-known that $\V P_{1}$, $\V P_{2}$,
$\V P_{3}$, and $\V P_{4}$ are all spectrally-equivalent if exact
solvers are used \cite{Elman_book}. Furthermore, if an exact Schur
complement inverse is employed, it can be shown for $\V P_{2}$, $\V P_{3}$,
and $\V P_{4}$ that any Krylov subspace iterative method with a Galerkin
property will require only a small number of iterations (two or three)
to converge to the exact answer \cite{Murphy1}.

As an alternative approximation to (\ref{sol_exact_p},\ref{sol_exact_u})
that is more accurate than the previous approximations, we consider
a fifth preconditioner closely-related to the Uzawa method, denoted
by $\V P_{5}$. The action of the inverse of this preconditioner $\V P_{5}^{-1}$
cannot easily be written in block-factorization form so we present
in the form of pseudo-code:
\begin{enumerate}
\item Solve for $\V x_{\ub}^{*}={\widetilde{\Ab}}^{-1}\V b_{\ub}$ using
multigrid with initial guess $\V 0$. 
\item Estimate pressure as $\V x_{p}\approx-\widetilde{\mathcal{S}}^{-1}(\Db\V x_{\ub}^{*}+\V b_{p})$. 
\item Estimate velocity as $\V x_{\ub}\approx\widetilde{\Ab}^{-1}(\V b_{\ub}-\Gb\V x_{p})$
using a multigrid solver, starting with $\V x_{\ub}^{*}$ as an initial
guess.
\end{enumerate}
If exact solvers are employed the only approximation made in $\V P_{5}$
is the approximation $\V S^{-1}\approx\widetilde{\mathcal{S}}^{-1}$,
and as such we expect it to be the best approximation to $\M M^{-1}$.
It is, however, also the most expensive of the five preconditioners
because it involves two applications of $\widetilde{\Ab}^{-1}$. Our
goal will be to investigate how well these preconditioners perform
in practice with inexact subsolvers.

\section{\label{sec:Staggered}Numerical Implementation}

In this section we specialize the relatively general preconditioners
from the previous section to a specific second-order conservative
finite-volume discretization of the time-dependent Stokes equations
on a uniform rectangular grid. We do not discuss here the inclusion
of advection in the full Navier-Stokes equations. Schemes that handle
advection explicitly using a non-dissipative spatial discretization
are described in detail in Refs. \cite{LLNS_Staggered,LowMachExplicit},
and Ref. \cite{NonProjection_Griffith} describes a particular higher-order
upwind scheme for uniform staggered-grids.

\subsection{Staggered-grid Discretization}

For our numerical investigations of the various preconditioners we
employ the well-known staggered-grid or MAC discretization of the
Stokes equations \cite{HarWel65,Guermond}. This is a conservative
discretization that is uniformly div-stable \cite{StaggeredGrid_Stability,Olshanskii}.
The scheme defines the degree of freedoms at staggered locations.
Specifically, scalar variables including pressure and density are
defined at cell centers, while components of vector variables including
velocity components are defined at the corresponding faces of the
grid \cite{NonProjection_Griffith,LLNS_Staggered}. For illustration,
we assume that the domain $\Omega$ is rectangular and there are $n_{x}$
cells along the $x$ direction and $n_{y}$ cells along the $y$ direction,
with periodic, no-slip (e.g., $\V u=0$ along a boundary) or free-slip
(e.g., $v=0$ and $\partial u/\partial y=0$ along the south boundary)
boundary conditions specified at each of the domain boundaries. For
simplicity, we further assume that the grid spacing along the different
directions is constant, $h_{x}=h_{y}=h$.

The divergence of ${\ub}=(u,v)^{T}$ is approximated at cell centers
by $\Db{\ub}=D^{x}u+D^{y}v$ with 
\[
(D^{x}u)_{i,j}=h^{-1}\left(u_{i+\myhalf,j}-u_{i-\myhalf,j}\right),\quad\quad(D^{y}v)_{i,j}=h^{-1}\left(v_{i,j+\myhalf}-v_{i,j-\myhalf}\right).
\]
The gradient of $p$ is approximated at the $x$ and $y$ edges of
the grid cells (faces in three dimensions) by $\Gb p=(G^{x}p,\, G^{y}p)^{T}$
with 
\[
(G^{x}p)_{i-\myhalf,j}=h^{-1}\left(p_{i,j}-p_{i-1,j}\right),\quad\quad(G^{y}p)_{i,j-\myhalf}=h^{-1}\left(p_{i,j}-p_{i,j-1}\right).
\]
For periodic domains or where a homogeneous Dirichlet condition is
specified for the normal component of velocity at physical boundaries,
the staggered discretization satisfies $\Db=-\Gb^{*}$. Note that
$\Db\Gb=\M L_{p}$, where $\M L_{p}$ is the standard (five-point
in two dimensions, seven-point in three dimensions) centered finite
difference Laplacian.

For constant viscosity, the finite difference approximation to the
vector Laplacian $\grad^{2}{\ub}$ is denoted as $\V L{\ub}=(L^{x}u,\, L^{y}v)$.
In the interior of the domain, $\grad^{2}u$ is discretized using
the standard five-point discrete Laplacian. In the presence of physical
boundaries, $L^{x}u$ is defined at all interior edges/faces where
$u$ are defined, and $L^{y}v$ is defined at all interior edges/faces
where $v$ are defined. The finite-difference stencils for tangential
velocities next to no-slip and free-slip boundaries are modified to
account for the boundary conditions, as described in Refs. \cite{NonProjection_Griffith,LLNS_Staggered}.
Note that for constant viscosity, if one uses the Laplacian form of
the viscous term, the different components of velocity are uncoupled.

When the viscosity is not a constant, the strain tensor form of the
viscous term is needed, for which 
\begin{equation}
\V L_{\V{\mu}}\V u=\grad\cdot\taub(\ub)={\displaystyle \left[\begin{array}{l}
2\frac{\partial}{\partial x}\left({\mu}\frac{\partial u}{\partial x}\right)+\frac{\partial}{\partial y}\left({\mu}\frac{\partial u}{\partial y}+{\mu}\frac{\partial v}{\partial x}\right)\\
2\frac{\partial}{\partial y}\left({\mu}\frac{\partial v}{\partial y}\right)+\frac{\partial}{\partial x}\left({\mu}\frac{\partial v}{\partial x}+{\mu}\frac{\partial u}{\partial y}\right)
\end{array}\right].}\label{div_stress}
\end{equation}
The discretization of $\grad\cdot\taub(\ub)$ is constructed using
standard (staggered) centered second-order differences to give the
discrete viscous operator $\V L_{\V{\mu}}$. Note that even for constant
viscosity, there is coupling between the velocity components in (\ref{div_stress}).
For the staggered discretization that we employ here, it can be shown
that for constant viscosity $\mu_{0}$, the viscous operator degenerates
to a Laplacian, $\V L_{\V{\mu}}\V u=\mu_{0}\V L\V u$, if $\Db{\ub}=0$.
That is, for constant viscosity incompressible flow the solution of
the Stokes system is not affected by the choice of the form of the
viscous term (By contrast, in fractional step methods, the unprojected
velocity and therefore the projected velocity is affected by this
choice). However, the Stokes solver is in general affected by the
choice of the viscous term, even for constant viscosity. As described
in Ref. \cite{LowMachExplicit}, centered differences for the viscous
fluxes that require values outside of the physical domain are replaced
by one-sided differences that only use values from the interior cell
bordering the boundary and boundary values. The tangential momentum
flux is set to zero for any faces of the corresponding control volume
that lie on a free-slip boundary, and values in cells outside of the
physical domain are never required. The overall discretization is
spatially globally second-order accurate.

We build the discrete velocity operator $\Ab=\theta{\V{\rho}}-\V L_{\V{\mu}}$
from the above centered finite-difference operators. We assume that
the density $\rho$ is specified at the cell centers. The density
matrix ${\V{\rho}}$ is constructed by defining the discrete momentum
density $\rho{\ub}$ at the cell faces, where the corresponding velocity
components are defined. Here we follow Ref. \cite{LowMachExplicit}
and average the density from cell centers to cell faces, 
\[
({\V{\rho}}{\ub})_{i+\myhalf,j}=\left(\frac{\rho_{i,j}+\rho_{i+1,j}}{2}\right)u_{i+\myhalf,j}\quad\mbox{and}\quad({\V{\rho}}{\ub})_{i,j+\myhalf}=\left(\frac{\rho_{i,j}+\rho_{i,j+1}}{2}\right)v_{i,j+\myhalf},
\]
giving a diagonal density matrix ${\V{\rho}}$ with the interpolated
face-centered densities along the diagonal. We will assume here that
the shear ${\V{\mu}}$ and bulk ${\V{\gamma}}$ viscosities are specified
at the cell centers; typically they are an explicit function of other
scalar variables such as density, temperature, and/or composition.
The matrices ${\V{\mu}}$ and ${\V{\gamma}}$ that appear in the approximation
to the Schur complement {[}e.g. Eq. (\ref{eq:Schur_app_bulk}){]}
are diagonal matrices containing the cell-centered values of the shear
and bulk viscosities. The discretization of the viscous operator $\V L_{\V{\mu}}$
requires a shear viscosity at both cell-centers and nodes (edges in
three dimensions). The value of ${\mu}$ at a node is set to be the
average of the four neighboring cell-centered values \cite{LowMachExplicit}.
Note that the Schur complement approximation uses only the cell-centered
and not the node-centered viscosities.

\subsection{Krylov Solver}

Having defined the discrete operators appearing in the Stokes system
(\ref{Saddle_sys}), we briefly discuss some issues that arise when
solving this saddle-point problem using an iterative Krylov solver.
The basic operation required by the Krylov solver is computing $\M M\V x$,
which amounts to a straightforward direct evaluation of the appropriate
finite-difference stencils at every interior face and every cell center
in the computational grid.

Application of any of the preconditioners requires implementing approximate
solvers for the pressure and velocity subproblems (i.e., application
of $\widetilde{\V L}_{\V{\rho}}^{-1}$and $\widetilde{\Ab}^{-1}$).
Here we employ geometric multigrid techniques to implement these solvers.
For the cell-centered pressure solver we use standard variable-coefficient
Poisson multigrid solvers \cite{almgren-iamr}. For the face-centered
velocity solver we develop a vector variant of a standard scalar Helmholtz
solver based on a generalized red-black Gauss-Seidel smoother. The
details of the multigrid algorithms are given in Appendix \ref{sec:Multigrid}.
Our implementation is based on the Fortran version of the BoxLib library
\cite{BoxLibWhitePaper}.

Note that because the preconditioners in Krylov methods are applied
to a residual-correction system, zero is an appropriate initial guess
for the multigrid subsolvers. Note also that for certain choices of
boundary conditions the pressure subproblem has a null space of constant
vectors. Similarly, for periodic steady-state problems the velocity
equation has a ($d$-dimensional) null-space of all constant velocity
fields. In these cases some care is needed in the implementation of
the preconditioners to ensure that the null-space is handled consistently
and the mean pressure and momentum are kept constant at certain prescribed
values (e.g., zero). When non-homogeneous velocity boundary conditions
are specified, the viscous operator is an affine operator because
the viscous stencils near the boundary use the specified boundary
values. Because Krylov solvers require a linear rather than an affine
operator, we apply the Krylov solver to the homogeneous form of the
Stokes problem by subtracting the boundary terms in a pre-processing
step.

Here we employ left preconditioning and apply the iterative solver
to the preconditioned system $\M P^{-1}\M M\V x=\M P^{-1}\V b$. The
convergence criterion for the Krylov solver is therefore most naturally
expressed in terms of either the absolute or relative reduction in
the magnitude of the preconditioned residual $r_{P}=\norm{\M P^{-1}\left(\M M\V x-\V b\right)}_{2}$.
A more robust alternative is to base convergence criteria on the value
of the unpreconditioned (true) residual $r=\norm{\M M\V x-\V b}_{2}$.
For the problems studies here we observe $r$ and $r_{P}$ to exhibit
similar convergence for well-scaled problems.

\subsection{Rescaling of the Linear System}

Another issue that arises when solving saddle-point problems is that
of scaling of the system to minimize the loss of floating point precision
occurs when adding the different terms. This is particularly important
when the equations are solved in dimensional form (with physical units),
but can also be important even when the equations are non-dimensionalized.
We consider re-scaling the velocity equations by some constant $c$
and re-scaling the pressure unknowns by the same factor, to obtain
the rescaled system
\begin{equation}
\left(\begin{array}{cc}
c\Ab & \Gb\\
-\Db & \V 0
\end{array}\right)\left(\begin{array}{c}
\V x_{\ub}\\
c\V x_{p}
\end{array}\right)=\left(\begin{array}{c}
c\V b_{\ub}\\
\V b_{p}
\end{array}\right).\label{eq:Stokes_rescaled}
\end{equation}
Intuitively, a well-scaled Stokes system is one in which both velocity-like
and pressure-like unknowns have elements of similar typical magnitude.
In order not to lose precision when evaluating the left-hand-side
of the velocity equation, we would like the viscous $c\Ab\V x_{\ub}$
and pressure $\Gb\V x_{p}$ contributions to have similar magnitude.
This suggests choosing $c\mu_{0}h^{-2}\sim h^{-1}$, giving $c\sim h/\mu_{0}$,
where $\mu_{0}$ is the typical magnitude of viscosity. Numerical
experiments confirm that rescaling the viscosity from a typical value
$\mu_{0}$ to $c\mu_{0}\sim h$ can dramatically improve the numerical
conditioning of the Stokes system. Note that no rescaling of the divergence
constraint is necessary since $\Db\V x_{\ub}$ has magnitude $\sim h^{-1}$
as the rest of the terms. Similarly, using equal weighting for velocity
and pressure residuals in the residual-space inner product will not
pose any problems because the two components of the residual $c\V b_{\ub}$
and $\V b_{p}$ already have similar magnitude. If there is a very
broad range of viscosities present in the problem, a uniform rescaling
of the equations will not be sufficient and diagonal scaling matrices
should be used to rescale the velocity and pressure separately, see
Eq. (31) in Ref. \cite{LargeViscosityJump} for a specific formulation.
To avoid loss of accuracy, in extreme cases extended precision arithmetic
may need to be used in the solver \cite{LargeViscosityJump}. 

In the numerical experiments reported in the next section we utilize
dimensionless well-scaled values ($h=1$, $\mu_{0}$=1, $\rho_{0}=1$)
for all of the coefficients, so that no explicit rescaling of the
unknowns or the equations is required. We have verified that after
the rescaling (\ref{eq:Stokes_rescaled}) we get similar results for
other choices of reference values for the viscosity and grid spacing.
We emphasize that if the re-scaling is not applied (i.e., $c=1$ is
used), multiple linear algebra issues arise when solving the Stokes
problem. For example, GMRES may not converge, or if it converges,
there may be a large difference between the preconditioned and unpreconditioned
residuals, and lastly, the computed solution may have a large error
due to ill-conditioning.

\section{\label{sec:Results}Results}

In this section we perform detailed numerical experiments to determine
the most robust and efficient preconditioner over a broad range of
parameters. Because the preconditioned system is not necessarily symmetric,
as a Krylov solver for the saddle-point system (\ref{Saddle_sys})
we use the left-preconditioned GMRES (Generalized Minimal Residual)
method with a fixed restart frequency $m$ \cite{Saad1,Saad2}. A
more robust and flexible method is FGMRES (Flexible GMRES). In particular,
FGMRES allows the use preconditioners that are not necessarily constant
linear operators (e.g., another Krylov solver or a variable number
of multigrid cycles). A notable drawback of FGMRES is that it requires
twice the storage of GMRES. Since one of our goals is to develop solvers
for large-scale calculations, we do not consider FGMRES here and use
the less memory-intensive GMRES method. GMRES requires the storage
of $m$ vectors like $\V x$. For a $d$-dimensional regular grid
with $N$ cells, the memory storage requirement is thus at least $(d+1)mN$
floating-point numbers since there are $d$ velocity degrees of freedom
(DOFs) and one pressure DOF per grid cell. It is therefore important
to explore the use of restarts to reduce the memory requirements of
the Krylov solver. For testing purposes, we generate a random $\V x$
and then compute $\V b=\M M\V x$, with zero velocity along all non-slip
boundaries. Similar convergence (not shown) is observed for other
right-hand sides or boundary velocities.

The multigrid algorithms used in the pressure and velocity subsolvers
iteratively apply V cycles, each of which consists of successive hierarchical
restriction (coarsening), smoothing, and prolongation \cite{Briggs}.
In our tests, we will use the same constant number $n$ of V cycles
in both the pressure and the velocity solvers. This ensures that the
preconditioners are constant linear operators and allows for the use
of the GMRES method. The velocity (vector) multigrid V cycle has a
cost very similar to $d$ independent pressure (scalar) V cycles.
Therefore, as a proxy for the CPU cost of a single application of
the preconditioner we will use the number of \emph{scalar} multigrid
cycles. The cost of the pressure subsolver (application of $\widetilde{\V L}_{\V{\rho}}^{-1}$)
is $n$ scalar V cycles, and the cost of the velocity subsolver (application
of $\widetilde{\Ab}^{-1}$) is $d\cdot n$ scalar V cycles. All preconditioners
require at least one velocity solve per application; however, they
differ in whether they require a pressure Poisson solve for steady
flow.

A fundamental ``easy'' test problem we employ is constant coefficient
steady Stokes flow in a periodic domain or a domain with no-slip condition
along all boundaries. As a more challenging variable-coefficient test
problem we use a \emph{bubble} test, in which we embed a sphere (disk
in two dimensions) of one fluid in another fluid with different viscosity
and density. The bubble is placed in the center of a cubic (square)
domain of length $n_{c}$ cells with no-slip boundaries along all
domain boundaries. For the bubble problem, the viscous stress is taken
to be $\taub(\ub)={\mu}\left[\grad\ub+(\grad\ub)^{T}\right]$, and
the diagonal elements of the viscosity matrix $\V{\mu}$ and the density
matrix $\V{\rho}$ at cell centers are generated from the spatially-dependent
functions $\mu\left(\V x\right)=\mu_{0}f(\V x;r_{\mu})$ and $\rho\left(\V x\right)=\rho_{0}f(\V x;r_{\rho})$
respectively, where 
\begin{equation}
f(\V x;r)=\frac{1}{2}\left(r+1\right)+\frac{1}{2}\left(r-1\right)\mbox{tanh}\left(\frac{\mbox{d}(\V x,\Gamma)}{\epsilon}\right)+0.1\V R.\label{varcoeff_def}
\end{equation}
Here $r_{\mu}$ and $r_{\rho}$ are the viscosity and density contrast
ratios, $\Gamma$ is the interface, a sphere of radius $L/4$ placed
at the center of a cube with side of length $L=n_{c}h$, $\mbox{d}(\V x,\Gamma)$
is the distance function to the interface, $\epsilon=h$ is a smoothing
width used to avoid discontinuous jumps in the coefficients, and $\V R$
is a random number uniformly distributed in $\left(0,1\right)$. Here
we focus on the case of no bulk viscous term; we have also done tests
including a bulk viscous term and observed similar behavior. Unless
otherwise indicated, the bubble test is steady-state ($\theta=0$),
and we use a contrast ratio $r_{\mu}=r_{\rho}=100$; we have observed
similar behavior with a ratio of $1000$, at least for sufficiently
smooth jumps (i.e., sufficiently large $\epsilon$). It is important
to note that the types of preconditioners we use here have been shown
to effective even with much larger viscosity constrasts ($\sim10^{6}$)
in the context of geophysical flow problems \cite{VariableViscosity_FEM,VariableViscosity_FEM2,LargeViscosityJump}.
For the target applications we have in mind, such as phase-field models
of fluid mixtures, viscosity contrasts of $10^{2}-10^{3}$ are more
relevant (for example, the viscosity ratio of water and air is only
55, while the density ratio is 1000).

\subsection{Multigrid Subsolvers}

Before comparing the different preconditioners, we optimize the key
parameters in the multigrid pressure and velocity approximate subsolvers,
specifically, the number of smoothing (relaxation) sweeps per V cycle
and the number of V cycles per application of the preconditioner.

\subsubsection{The effect of the number of smoothing sweeps}

One of the key aspects of geometric multigrid is the smoother used
to perform relaxation of the error at each level of the multigrid
hierarchy. As explained in more detail in Appendix \ref{sec:Multigrid},
we employ a red-black Gauss-Seidel smoother. This ensures that all
components of the error are damped to some extent for constant-coefficient
problems, and, more importantly, makes the smoother highly parallelizable.
The optimal number of smoothing (relaxation) sweeps to be performed
at each multigrid level (we use the same number of sweeps going down
and up the multigrid hierarchy) has to be determined by numerical
experimentation.

\begin{figure}[!t]
\begin{centering}
\includegraphics[width=0.49\textwidth]{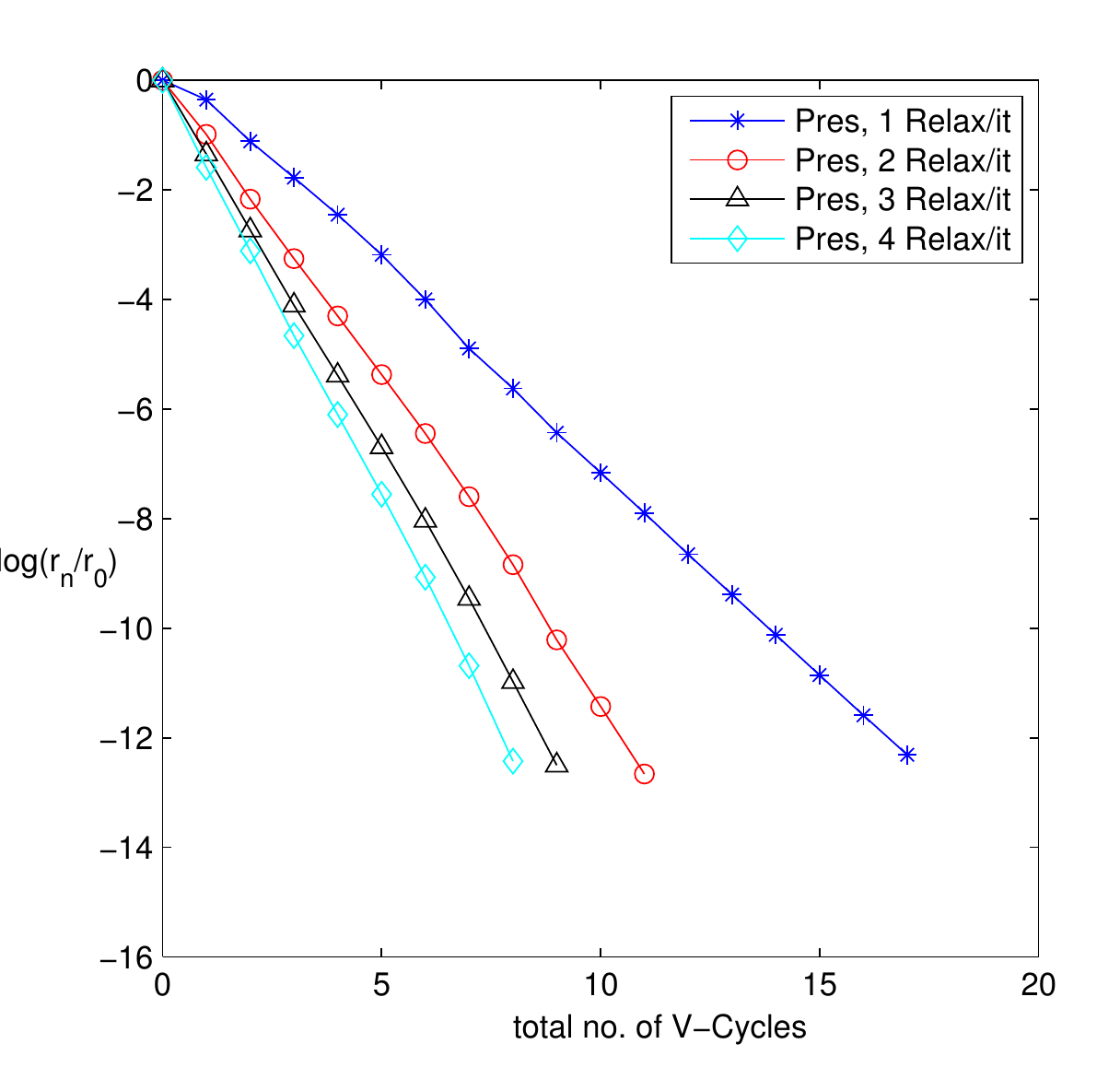}\includegraphics[width=0.49\textwidth]{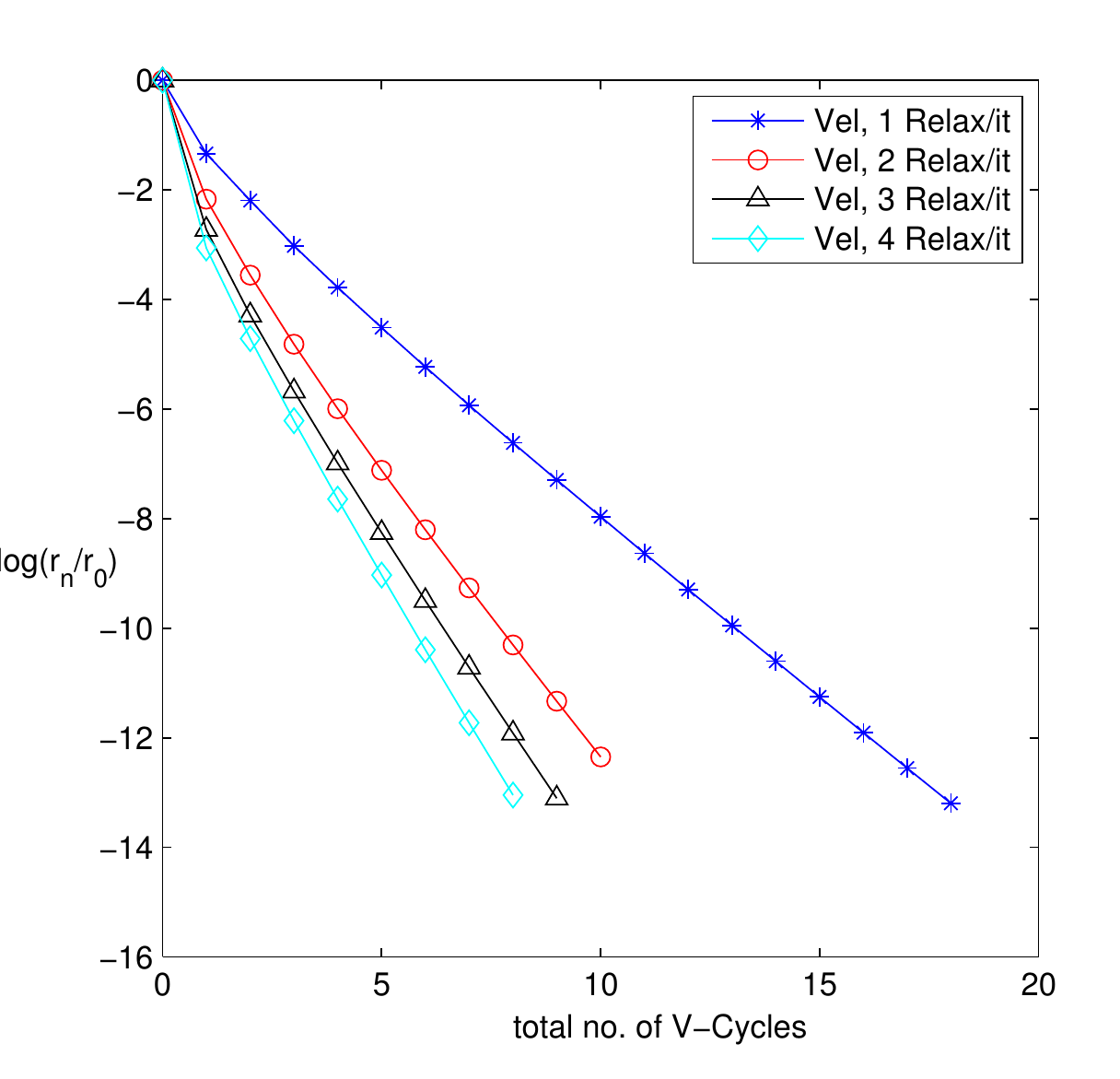}
\par\end{centering}

\begin{centering}
\includegraphics[width=0.49\textwidth]{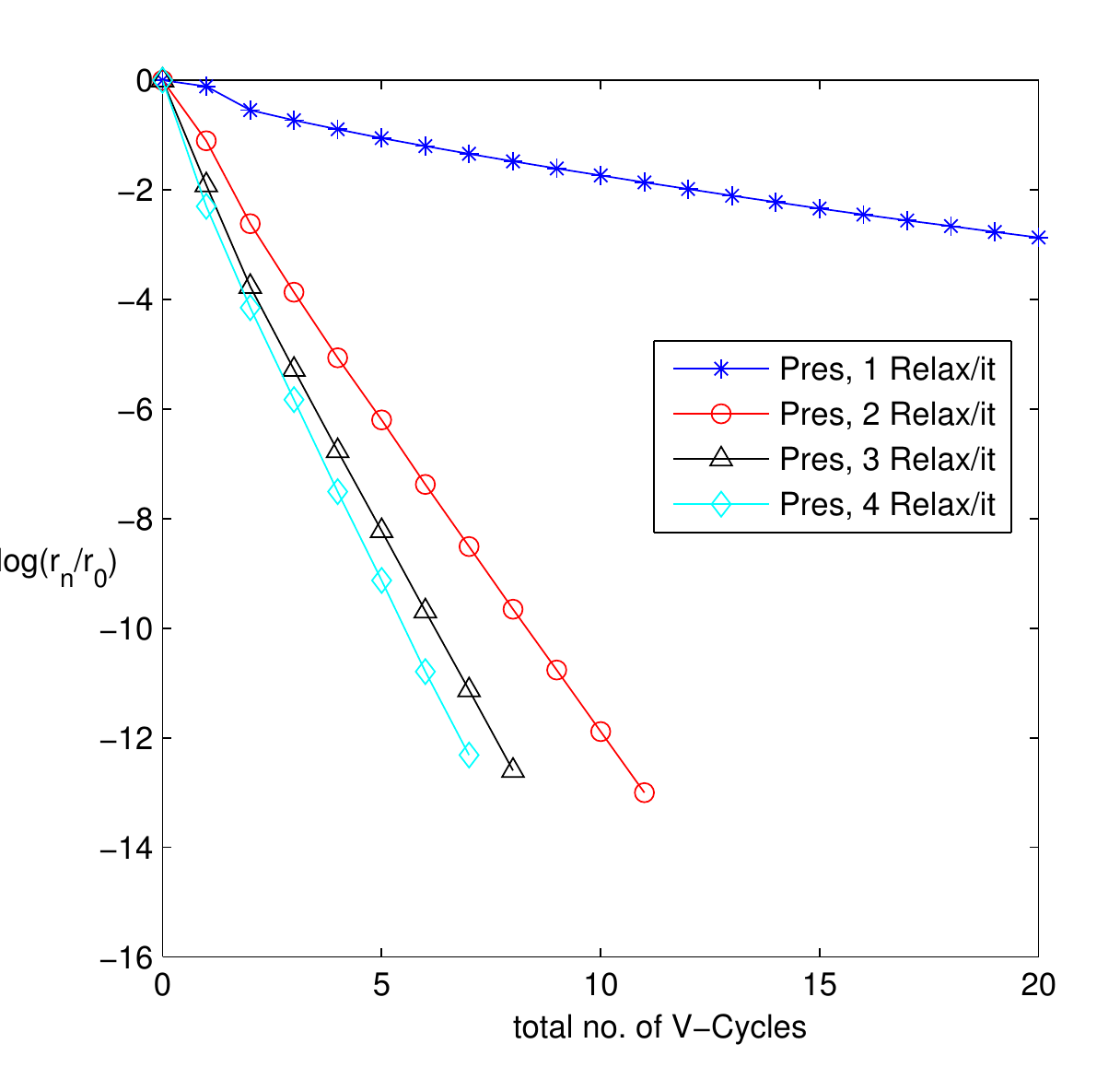}\includegraphics[width=0.49\textwidth]{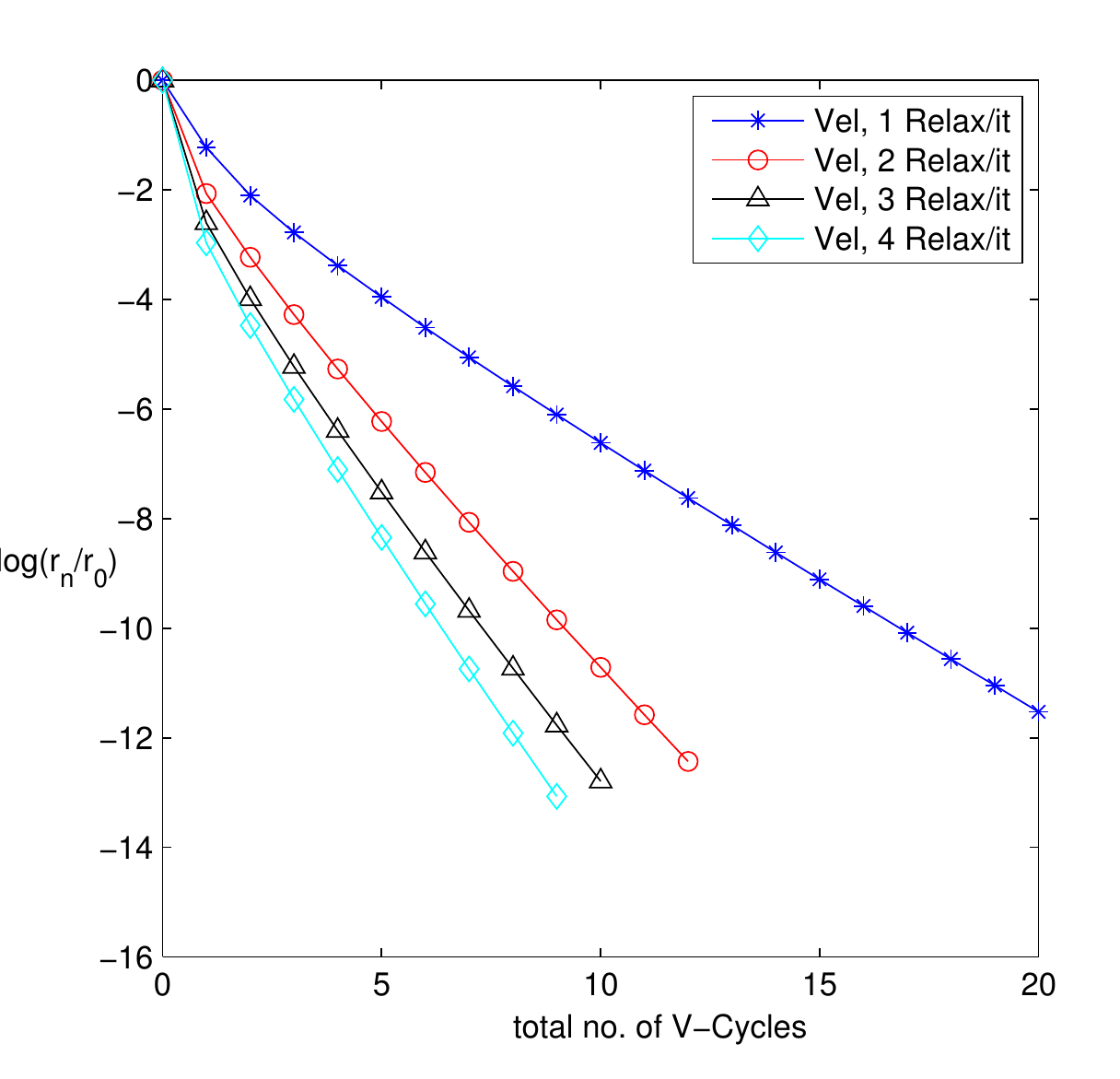}
\par\end{centering}

\caption{\label{fig:MGRelaxation}The log of the relative residual for the
pressure (left) and velocity (right) multigrid solvers as a function
of the number of multigrid V cycles, for different numbers of smoothing
(relaxation) sweeps. A constant coefficient steady Stokes problem
is solved on a $512^{2}$ grid in two dimensions (top panels), and
$128^{3}$ grid in three dimensions (bottom panels), with no-slip
conditions at all domain boundaries.}
\end{figure}

In Fig. \ref{fig:MGRelaxation} we show the convergence of the pressure
(left panels) and velocity multigrid solvers (right panels) for constant
viscosity but for the stress-tensor form of the viscous term (\ref{div_stress}).
In the upper row we show results in two dimensions, and in the lower
row we show results for three dimensions. Similar results are obtained
for different types of boundary conditions. We see a large increase
in the rate of convergence when increasing the number of smoothing
sweeps from one to two, and only a modest increase thereafter. Since
the cost of geometric multigrid is in large part dominated by smoother,
henceforth we use two applications of the smoother at each level of
the multigrid hierarchy in each V cycle.

The speed of convergence of the multigrid iteration for the component
solvers, is the standard against which one should measure convergence
of the Krylov solver for the Stokes problem. As we can see in Fig.
\ref{fig:MGRelaxation}, each V cycle reduces the residual by at least
an order of magnitude, so that only about a dozen V cycles are needed
to reduce the residual to near roundoff. Therefore, a Stokes solver
that uses only $10(d+1)$ scalar multigrid cycles to reduce the residual
by more then 10 orders of magnitude should be considered excellent.

\subsubsection{The effect of the number of multigrid cycles}

As detailed in Appendix \ref{sec:AnalysisEigenvalues}, for constant-coefficient
Stokes problems with periodic boundaries, GMRES will converge in a
single iteration with preconditioner $\V P_{1}$ and in two iterations
with $\V P_{2}$ and $\V P_{3}$. The same holds for any choice of
boundary condition for the time-dependent case in the inviscid limit.
However, in the majority of cases of interest, multiple GMRES iterations
will be required even if the subsolvers are exact. It is therefore
important to explore the use of inexact pressure and velocity solvers.
Specifically, it is important to determine the optimal number of multigrid
V cycles per application of the preconditioner. We use the preconditioner
$\V P_{1}$ for these tests; similar results are observed for all
of the preconditioners.

\begin{figure}[!t]
\begin{centering}
\includegraphics[width=0.49\textwidth]{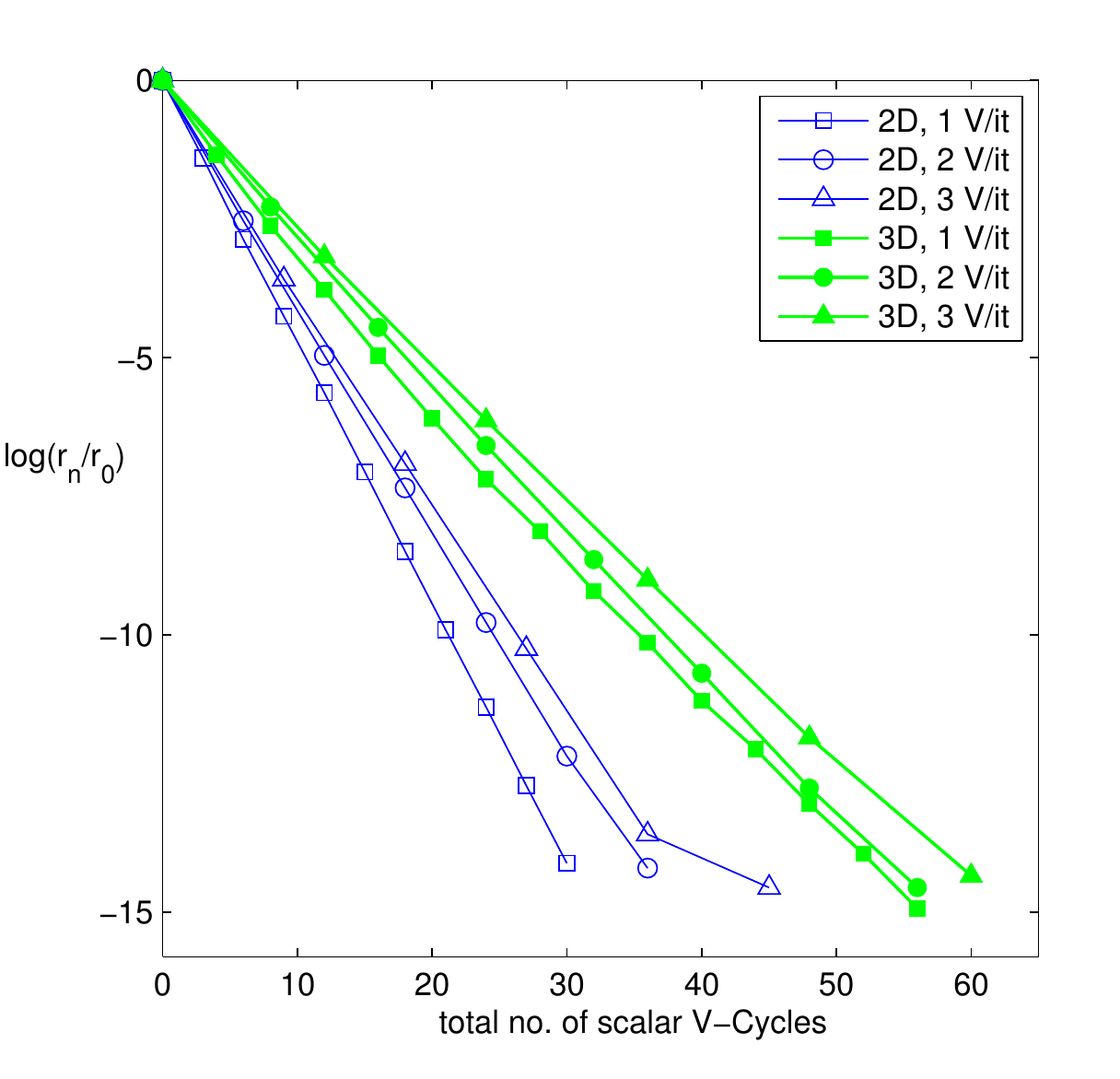}\includegraphics[width=0.49\textwidth]{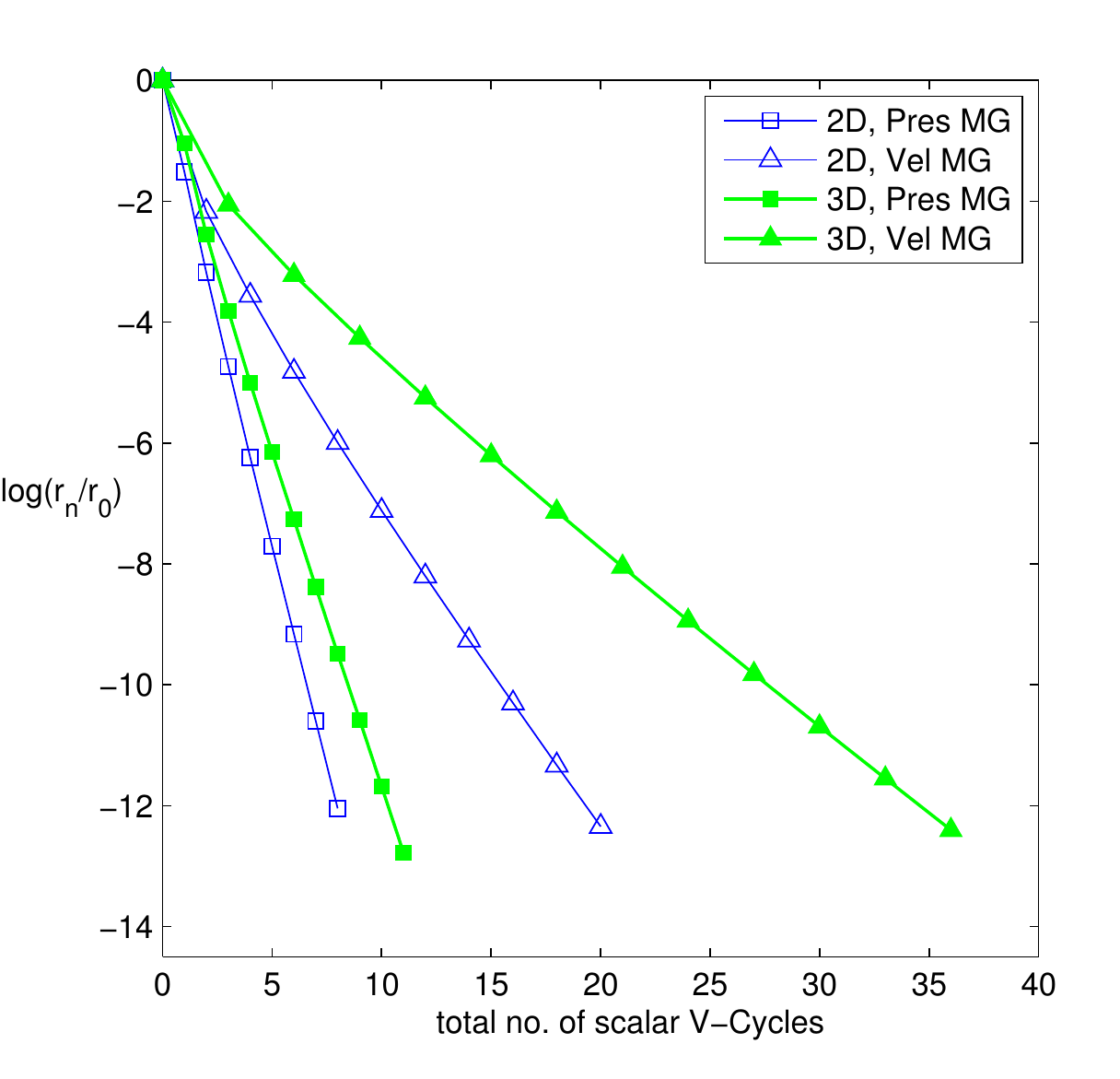}
\par\end{centering}

\begin{centering}
\includegraphics[width=0.49\textwidth]{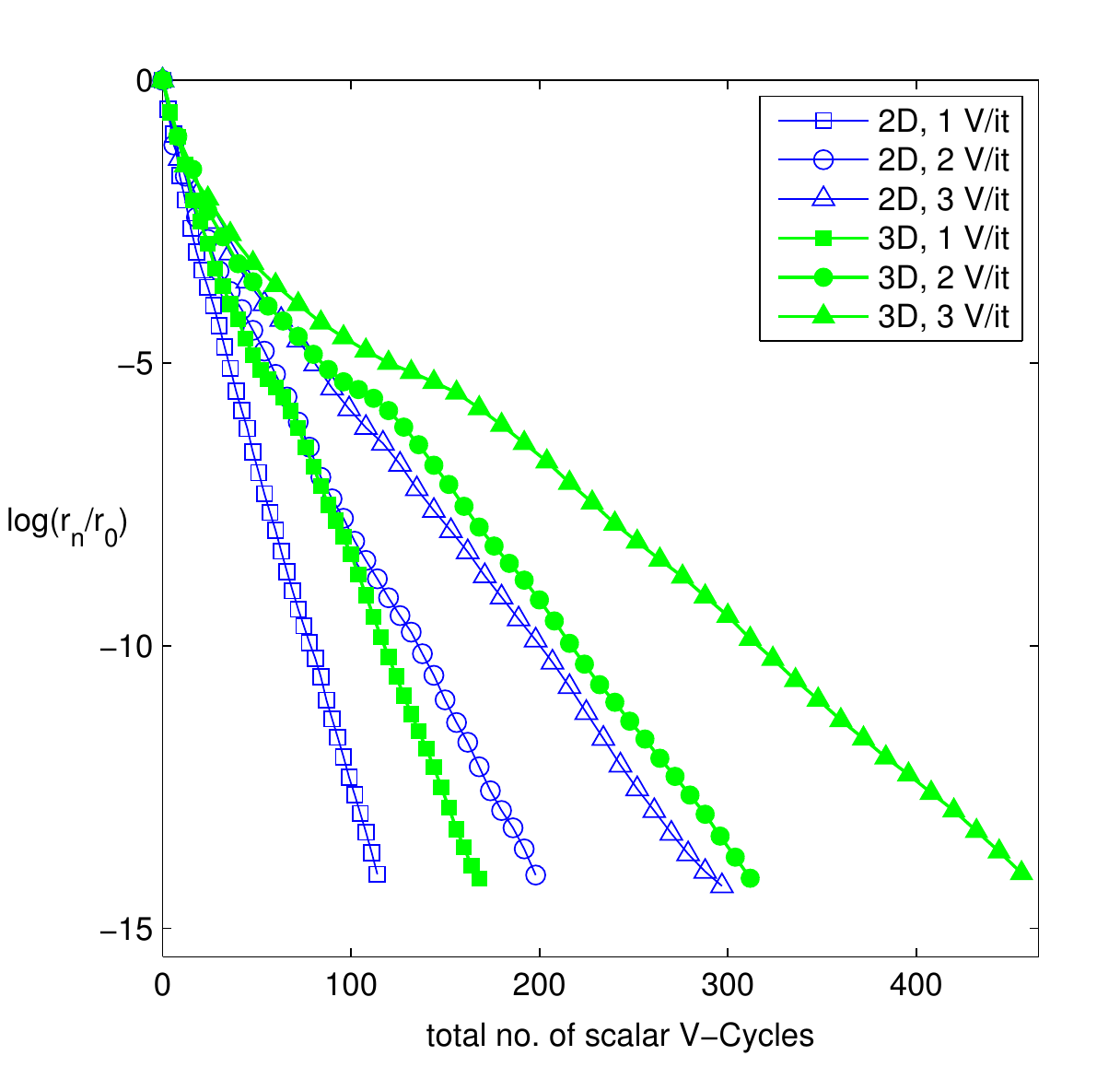}\includegraphics[width=0.49\textwidth]{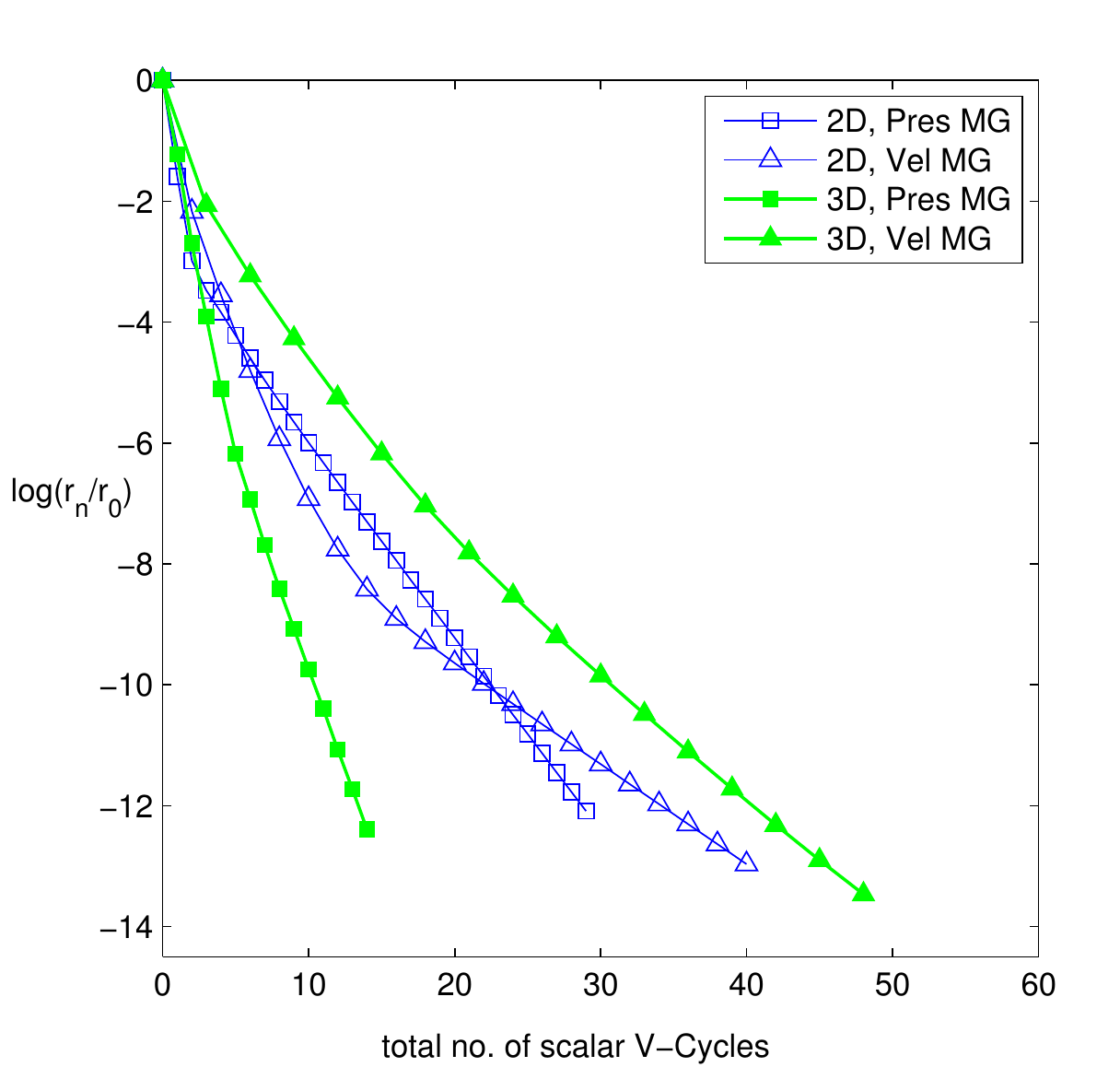}
\par\end{centering}

\caption{\label{fig:MGVCycles}The relative residual (on a log scale) as a
function of the total number of scalar multigrid V cycles, for different
number of multigrid cycles per application of the pressure and velocity
subsolvers in preconditioner $\V P_{1}$. GMRES convergence is shown
in the left panels, and pressure (squares) and velocity (triangles)
multigrid convergence is shown in the right panels, in both two ($512^{2}$
grid, empty symbols) and three ($128^{3}$ grid, filled symbols) dimensions.
Restarts are not employed in the GMRES solver. The top panels show
results for a constant-coefficient periodic steady-state Stokes problem,
and the bottom panels show results for the bubble test problem.}
\end{figure}

In the left panels of Fig. \ref{fig:MGVCycles} we show the convergence
of the relative preconditioned residual, as estimated by the GMRES
algorithm, for steady Stokes problems in two and three dimensions,
as a function of the total number of scalar V cycles. We recall that
the number of V cycles is a good proxy for the total computational
effort, so that the most rapid convergence in these plots corresponds
to the most efficient solvers. In the top left panel we show results
for constant viscosity but for the stress-tensor form of the viscous
term (\ref{div_stress}) for a periodic system, and in the bottom
left panel we show results for the variable-viscosity bubble problem
described earlier. In the corresponding right panels we show the convergence
of the pressure and velocity multigrid subsolvers on the same problem,
to serve as a reference point for what one may expect for a projection-like
split solver.

The top left panel in Fig. \ref{fig:MGVCycles} shows that for periodic
constant coefficient problems there is no significant difference between
using an exact subsolver (in effect, many V cycles per application
of the preconditioner), and using only a single V cycle in the preconditioner
but doing more GMRES iterations. This is not unexpected because the
standard multigrid algorithm takes the form of a simple Richardson
iterative solver, and we expect GMRES to perform at least as well
as Richardson iteration. Note that for more difficult Poisson problems,
such as problems with large jumps in the coefficients, it is well-known
that a Krylov solver preconditioned with multigrid is more robust
than multigrid alone, see for example the discussion in Ref. \cite{LargeViscosityJump}.

In the bottom left panel in Fig. \ref{fig:MGVCycles} we show the
convergence of GMRES for the variable-coefficient bubble problem,
which is typical of the behavior we observe when there are non-periodic
boundary conditions or variable coefficients. Similar behavior is
observed for the other preconditioners (not shown). The results clearly
demonstrate that when using exact subsolvers in the preconditioner
does not yield an exact solver for the Stokes problem, the extra cost
of performing more than a single V cycle of multigrid does not pay
off in terms of overall efficiency. The optimal rate of convergence
is observed when using only a single V cycle in the preconditioner.
We have observed no advantage to using a different number of cycles
in the pressure and velocity solvers, except for nearly inviscid problems
where performing more than one pressure cycle may be somewhat beneficial.
By comparing the left and right panels, we see that when using a single
multigrid cycle in the preconditioner the total number of scalar V
cycles is at most 2-3 times larger than that used in fractional step
(projection) methods (for example, $\sim50+15=65$ in three dimensions
for projection methods as seen in the right panel, and $\sim170$
cycles for coupled solver as seen in the left panel).

Based on these observations, henceforth we use only a single multigrid
cycle in the subsolvers employed by the preconditioners.

\subsection{Comparison of Preconditioners}

Having empirically determined the optimal settings for the pressure
and velocity subsolvers, we now turn to exploring the performance
of the different preconditioners. We begin by settling an issue regarding
the proper choice of sign in the upper/lower triangular and block-diagonal
preconditioners.

\subsubsection{The effect of the sign of Schur complement}

In the literature \cite{Ipsen,Murphy1}, the following Schur complement
based preconditioners have been proposed and studied, 
\[
\V P_{\pm}=\left(\begin{array}{cc}
\Ab & \V 0\\
-\Db & \pm(-\V D\Ab^{-1}\Gb)
\end{array}\right),
\]
where the sign of the lower diagonal block can be either positive
or negative. It was proven that $\V T_{+}=\V P_{+}^{-1}\V M$ satisfies
$(\V T_{+}-\Ib)(\V T_{+}+\Ib)=\V 0$ and $\V T_{-}=\V P_{-}^{-1}\V M$
satisfies $(\V T_{-}-\Ib)^{2}=\V 0$ \cite{Ipsen,Murphy1}. If a block-diagonal
preconditioner such as $\M P_{4}$ is used, by changing the sign it
is possible to either make the preconditioned operator symmetric but
indefinite (allowing the use of methods such as MINRES), or non-symmetric
but positive semi-definite \cite{FischerWathen}. 

Because the GMRES method possesses a Galerkin property \cite{Eiermann},
the total number of GMRES iterations is equal to the degree of the
characteristic polynomials of the preconditioned systems. Therefore,
using both $\V P_{+}^{-1}\V M$ and $\V P_{-}^{-1}\V M$, GMRES converges
in 2 iterations if the inverses of $\Ab$ and the Schur complement
are exact \cite{Murphy1}. However, when inexact subsolvers are employed,
we observe significant difference between the two choices of the sign
of the Schur complement. Our numerical results (not shown) indicate
that the preconditioners with \textquotedbl{}-\textquotedbl{} sign
in front of Schur complement give almost twice faster GMRES convergence
than those with \textquotedbl{}+\textquotedbl{} sign. This is consistent
with our original choice in Eqs. (\ref{eq:P2_S}) and (\ref{eq:P3_S}).

\subsubsection{The effect of restarts}

\begin{figure}[h]
\begin{centering}
\includegraphics[width=0.49\textwidth]{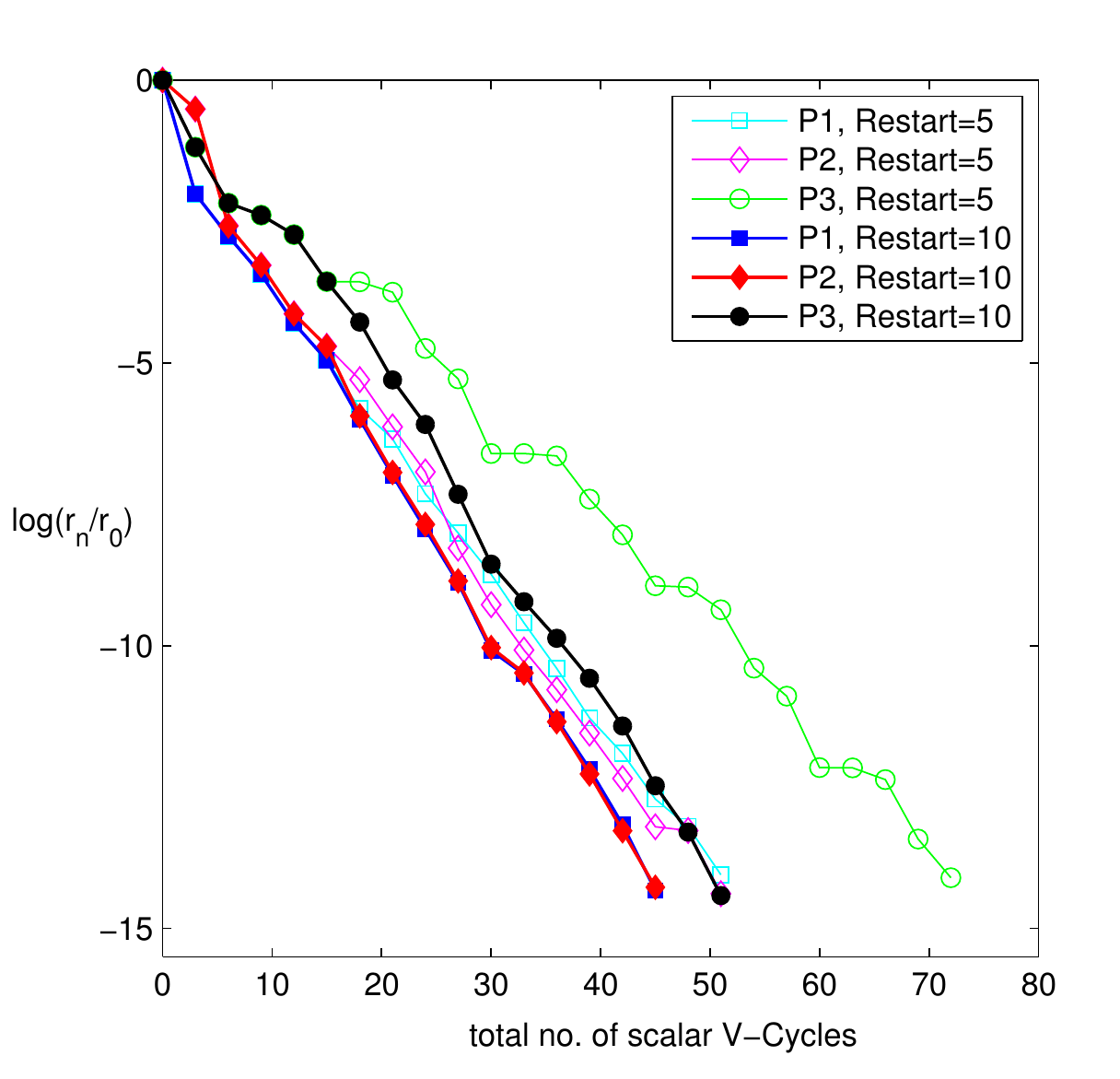}\includegraphics[width=0.49\textwidth]{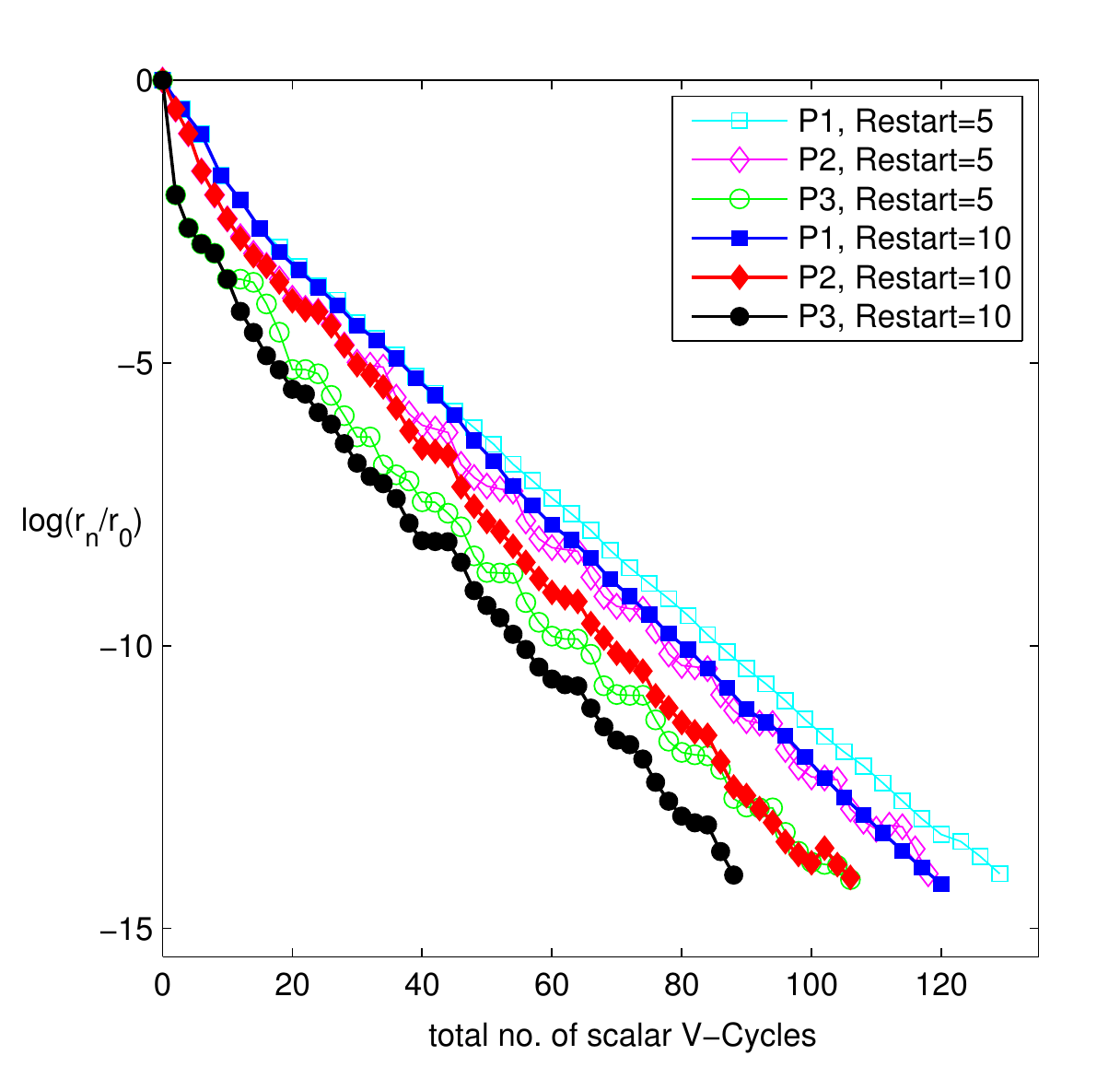}
\par\end{centering}

\caption{\label{fig:Restarts}The GMRES convergence history for preconditioners
$\V P_{1}$, $\V P_{2}$ and $\V P_{3}$ in two dimensions ($512^{2}$
grid) for the bubble test problem, for GMRES restart frequency 5 (open
symbols) and 10 (filled symbols). In the left panel we set the viscosity
to zero (unsteady inviscid flow) and in the right panel we set density
to zero (steady viscous flow).}
\end{figure}

For large-scale problems, particularly in three dimensions, the memory
requirements of the GMRES algorithm can be excessive. Restarts of
the GMRES iteration offer a simple way not only to avoid convergence
stalls, but also to limit the memory use. In Fig. \ref{fig:Restarts}
we compare the behavior of $\V P_{1}$, $\V P_{2}$ and $\V P_{3}$
for restart intervals of 5 or 10 GMRES iterations. In the left panel
of the figure we show the behavior for an inviscid time-dependent
bubble test problem ($\M L_{\V{\mu}}=\V 0$, relevant to simulations
of large Reynolds number flows) and in the right panel we show the
behavior for a steady Stokes bubble problem (relevant to small Reynolds
number flows). A two dimensional calculation is shown in the figure
but similar results are observed in three dimensions as well. In the
left panel of Fig. \ref{fig:Restarts} we see that the performance
of $\M P_{3}$ significantly deteriorates for the small restart interval
for the inviscid problem. In the right panel of the figure we see
some deterioration of the convergence for the small restart interval
for all three preconditioners.

Based on these results, henceforth we use a restart interval of 10
iterations.

\subsubsection{Comparisons of different preconditioners}

\begin{figure}[h]
\begin{centering}
\includegraphics[width=0.49\textwidth]{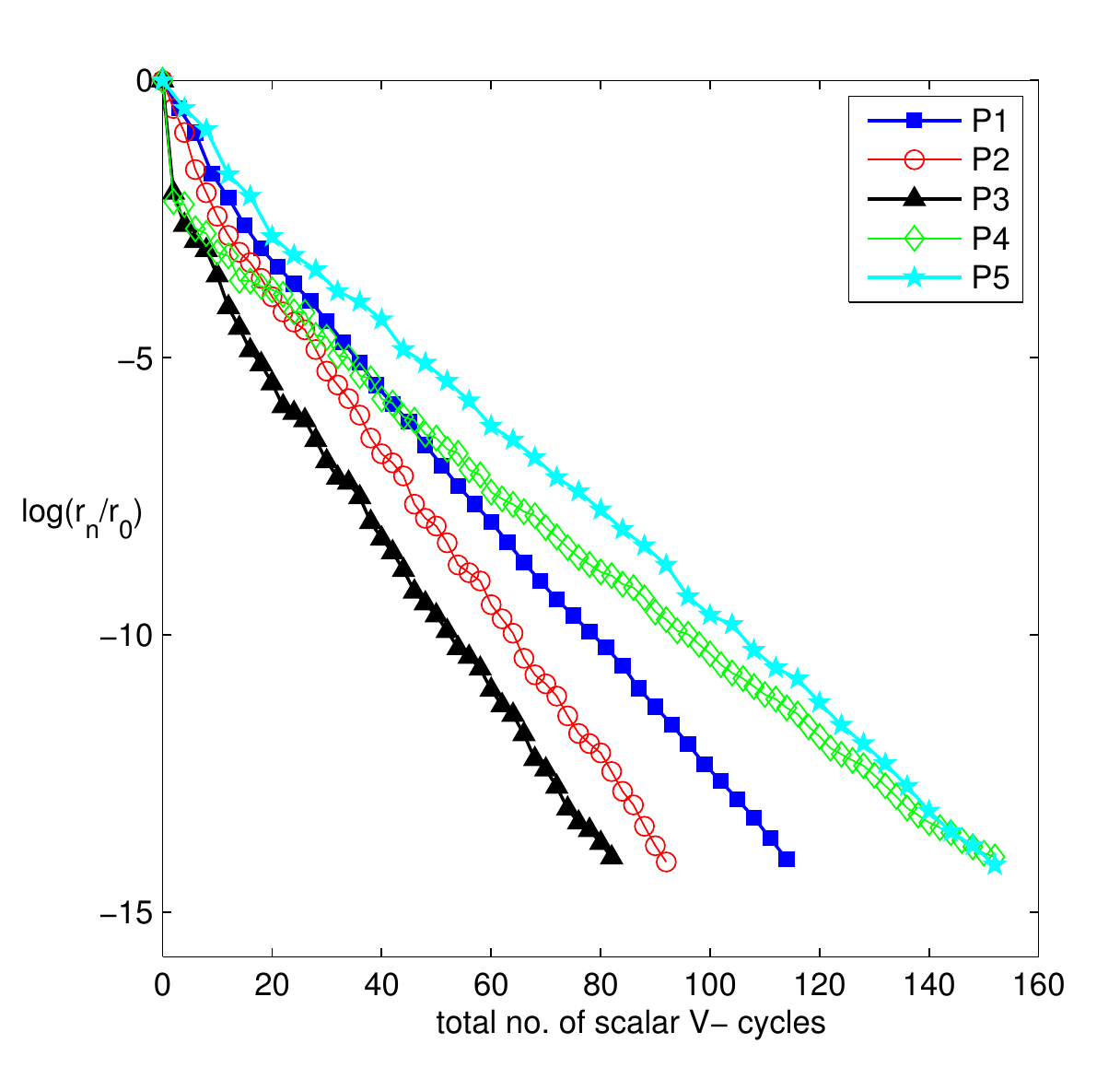}\includegraphics[width=0.49\textwidth]{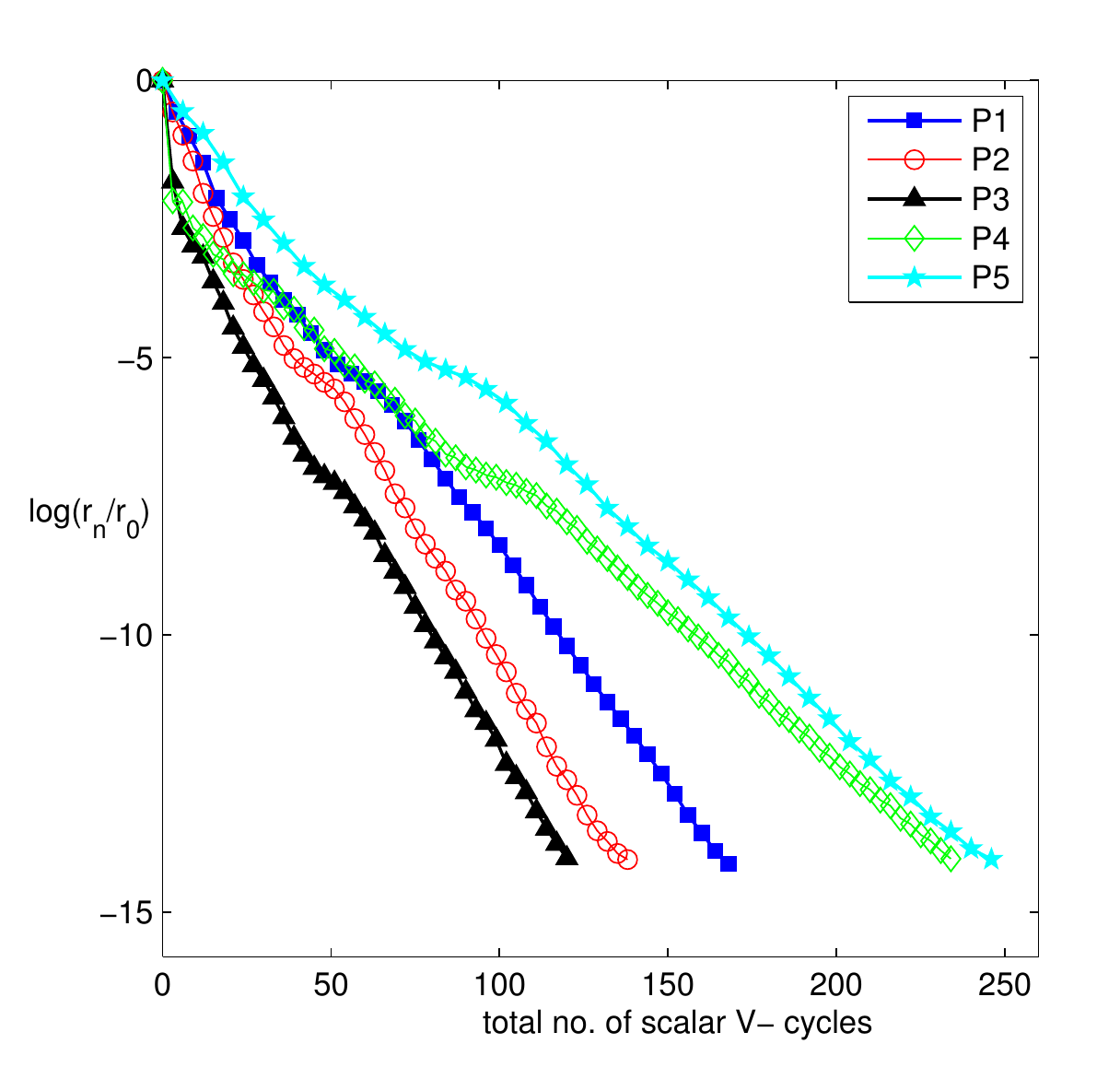}
\par\end{centering}

\caption{\label{fig:Comparison}The GMRES convergence history for preconditioners
$\V P_{1}$, $\V P_{2}$, $\V P_{3}$, $\V P_{4}$ and $\V P_{5}$
in two (left panel, $512^{2}$ grid) and three (right panel, $128^{3}$
grid) dimensions, for the bubble test problem. The restart frequency
is 10 GMRES iterations.}
\end{figure}

Having empirically determined the optimal subsolver settings and the
optimal sign of the Schur complement in the lower diagonal block of
the preconditioners, we can now compare the performance of the five
preconditioners on the bubble test problem in two and three dimensions.
The GMRES convergence results shown in Fig. \ref{fig:Comparison}
demonstrate that for steady Stokes problem the lower and upper triangular
preconditioners $\V P_{2}$ and $\V P_{3}$ yield the most efficient
GMRES solver. The projection preconditioner $\V P_{1}$ is seen to
give robust convergence but is less efficient for the steady flow
case because it requires one additional scalar (pressure) V cycle
per GMRES iteration. The results in the figure also clearly show that
$\V P_{4}$ and $\V P_{5}$ are much less efficient. This shows that
including an upper or lower triangular block in the Schur complement
based block preconditioners improves convergence, and also shows that
the extra work in $\V P_{5}$ over $\V P_{1}$ is not justified in
terms of overall efficiency, similarly to how the additional pressure
solve in $\V P_{1}$ does not yield improvement.

Based on these observations, henceforth we do not consider $\V P_{4}$
and $\V P_{5}$. Since we find very similar behavior between $\V P_{2}$
and $\V P_{3}$, while $\V P_{3}$ shows poorer behavior with frequent
restarts, henceforth we focus on examining in more detail the performance
of $\V P_{1}$ and $\V P_{2}$.

\subsection{Robustness}

In this section we examine in more detail the robustness of $\V P_{1}$
and $\V P_{2}$ under varying importance of the viscous contribution
to $\M A$, and changing problem size.

\subsubsection{Effects of viscous CFL number}

\begin{figure}[h]
\begin{centering}
\includegraphics[width=0.49\textwidth]{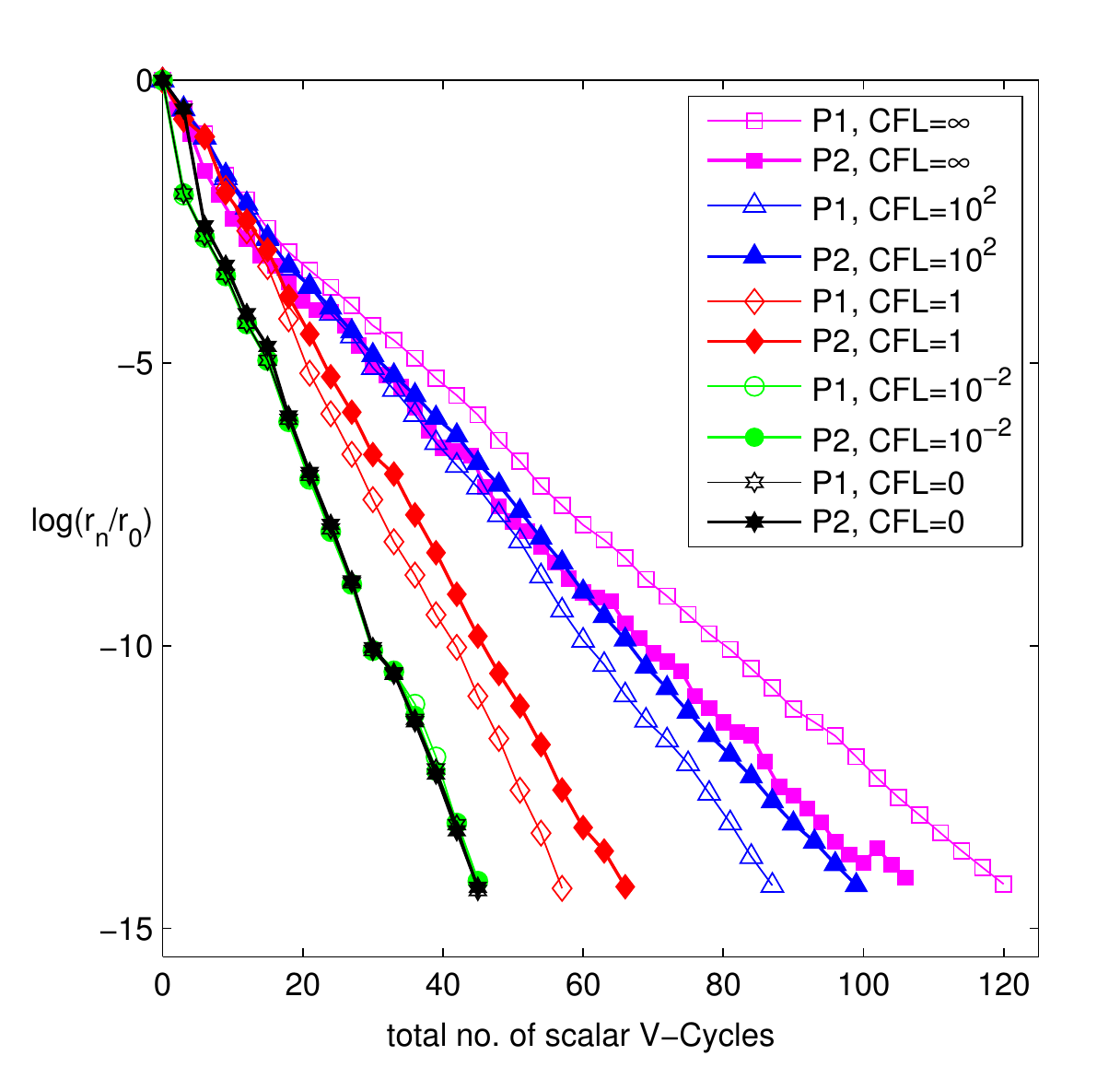}\includegraphics[width=0.49\textwidth]{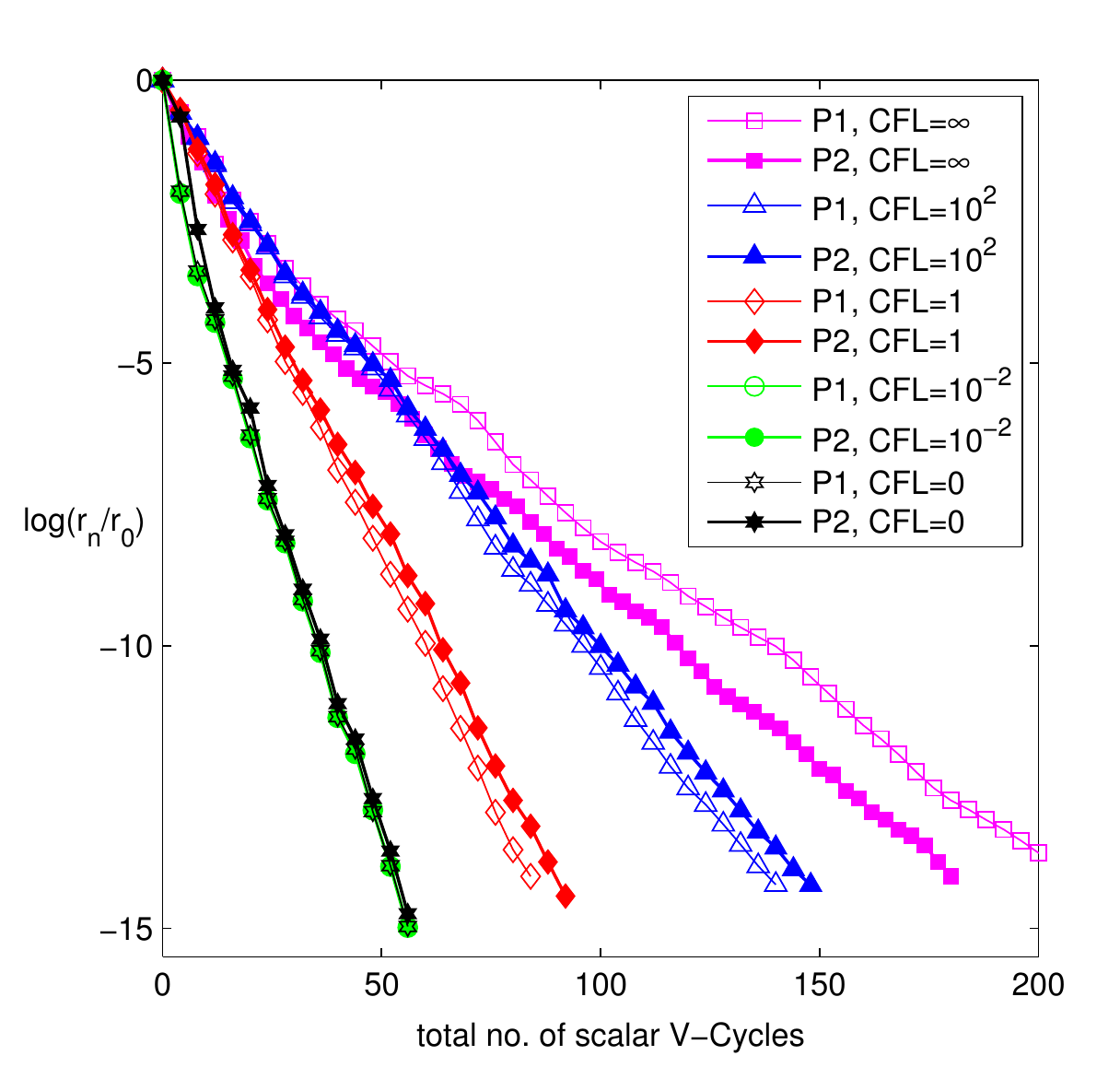}
\par\end{centering}

\caption{\label{fig:CFL} The GMRES convergence history for preconditioners
$\V P_{1}$ (empty symbols) and $\V P_{2}$ (filled symbols) in two
dimensions (left panel, $512^{2}$ grid) and in three dimensions (right
panel, $128^{3}$ grid) for the bubble test problem. We vary $\theta$
to change the viscous CFL number $\beta$ from the inviscid limit
$\beta=0$ to the steady limit $\beta\rightarrow\infty$.}
\end{figure}

One of the goals of our study is to design preconditioners that work
not just in the steady state limit but also for time-dependent problems.
While one can use a suitably-defined Reynolds number to measure the
importance of the inertial term $\theta{\V{\rho}}$ in $\M A$ relative
to the viscous term $\V L_{\V{\mu}}$, the best dimensionless number
to use for this is the viscous CFL number
\[
\beta=\frac{\nu_{0}}{\theta h^{2}}=\frac{\mu_{0}}{\theta\rho_{0}h^{2}}.
\]
A small $\beta\ll1$ indicates an easier problem where inertial effects
dominate, with $\beta=0$ corresponding to inviscid flow. A large
$\beta>n_{c}^{2}$ indicates a viscous-dominated problem, where $n_{c}$
is the grid size, with the hardest case being a steady-state problem
$\beta\rightarrow\infty$. In Fig. \ref{fig:CFL} we study the performance
of the GMRES Stokes solver for varying viscous CFL numbers for the
bubble test problem, in both two and three dimensions, for both preconditioners
$\V P_{1}$ and $\V P_{2}$. As expected, we see most rapid convergence
for $\beta=0$, and slowest convergence for $\beta\rightarrow\infty$.
For the steady state case $\theta=0$, we do not need a pressure Poisson
solve for $\V P_{2}$ and therefore this preconditioner is somewhat
more efficient than $\V P_{1}$. For intermediate $\beta$'s we get
somewhat better convergence for $\V P_{1}$, although the difference
is small. In our experience both preconditioners show rather robust
behavior for varying viscous CFL number, viscosity and density contrast
ratios, and different combinations of boundary conditions (periodic,
free-slip, or no-slip).

\subsubsection{Effects of problem size}

\begin{figure}[!th]
\begin{centering}
\includegraphics[width=0.49\textwidth]{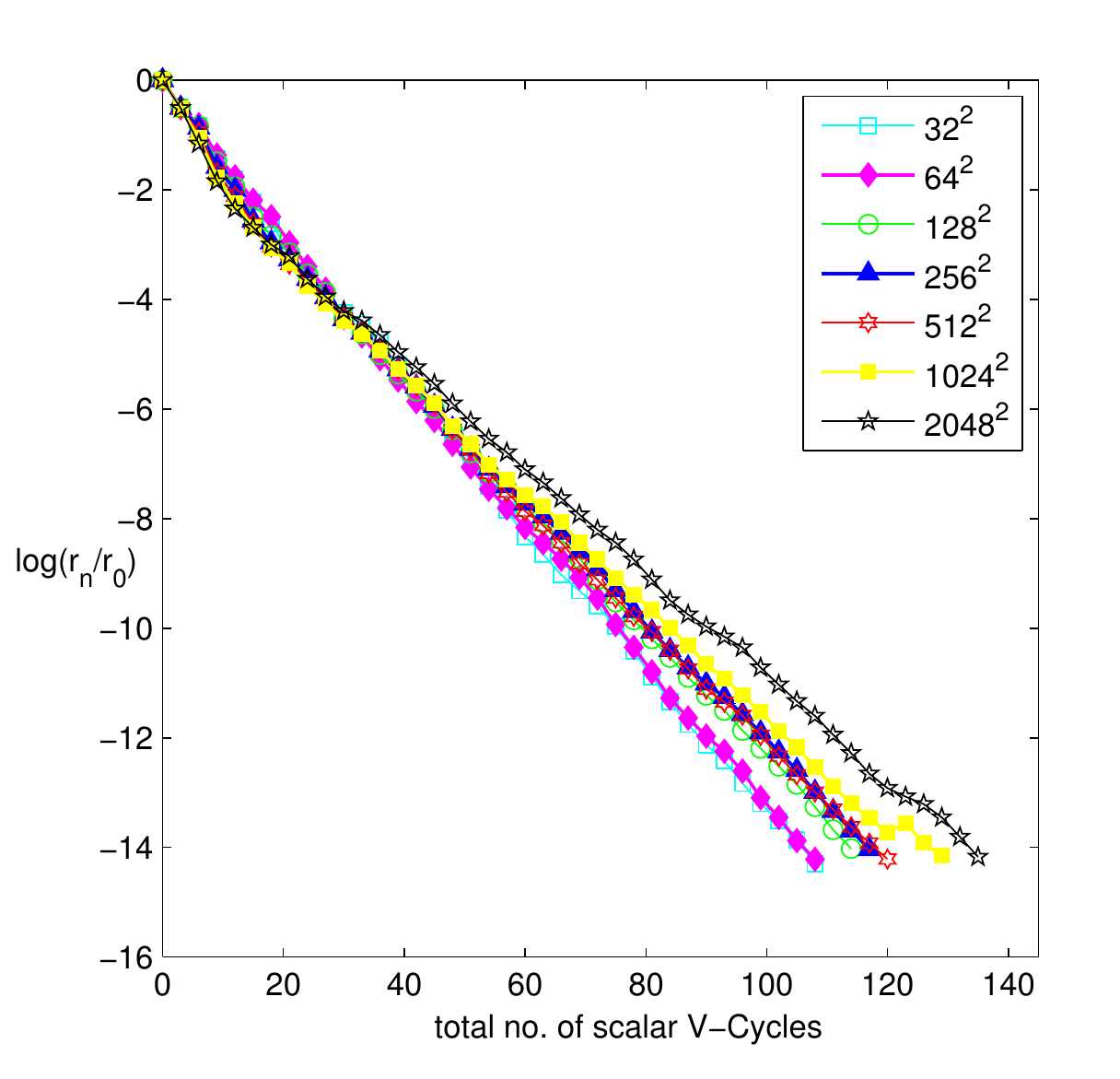}\includegraphics[width=0.49\textwidth]{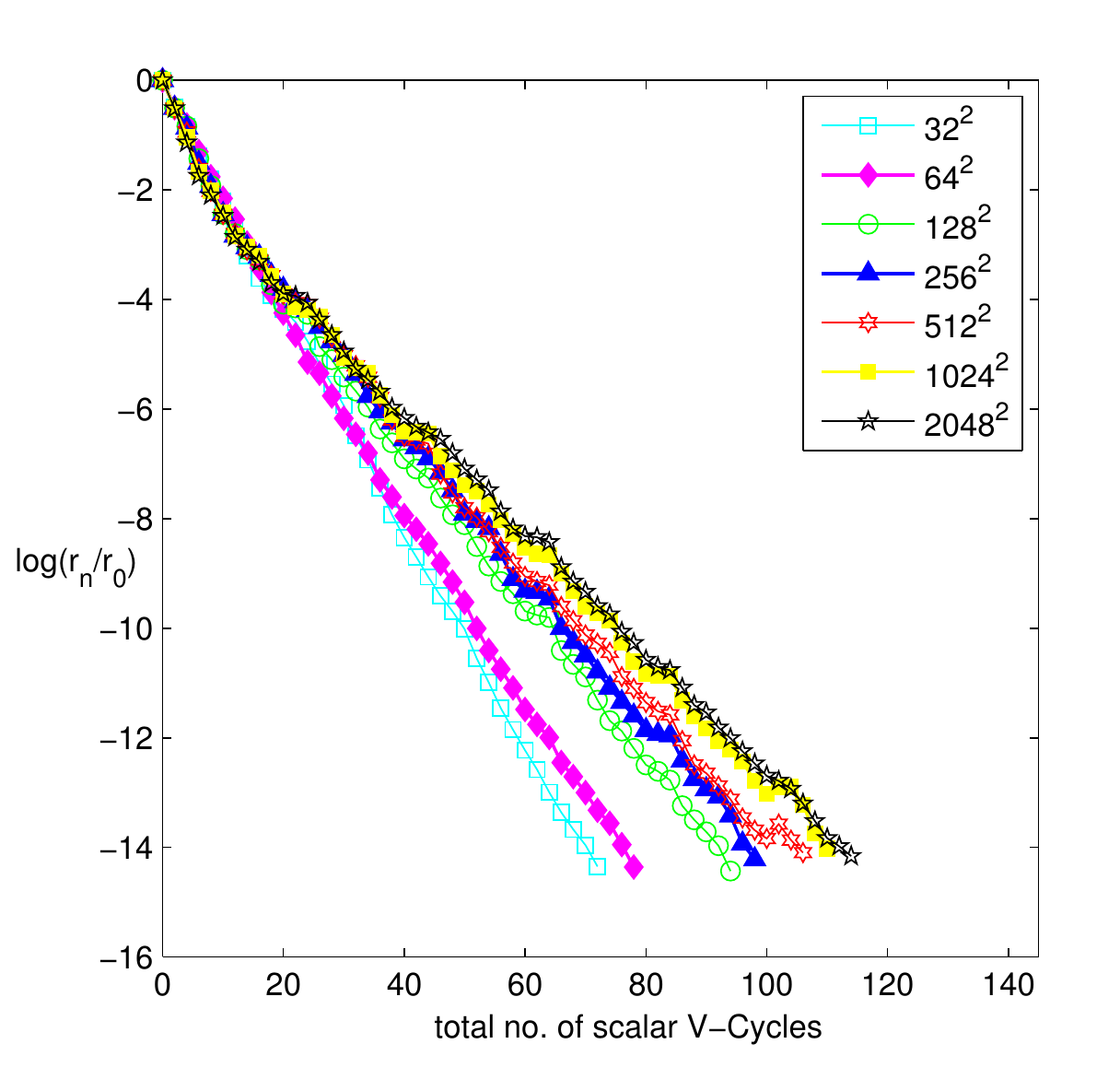}
\par\end{centering}

\begin{centering}
\includegraphics[width=0.49\textwidth]{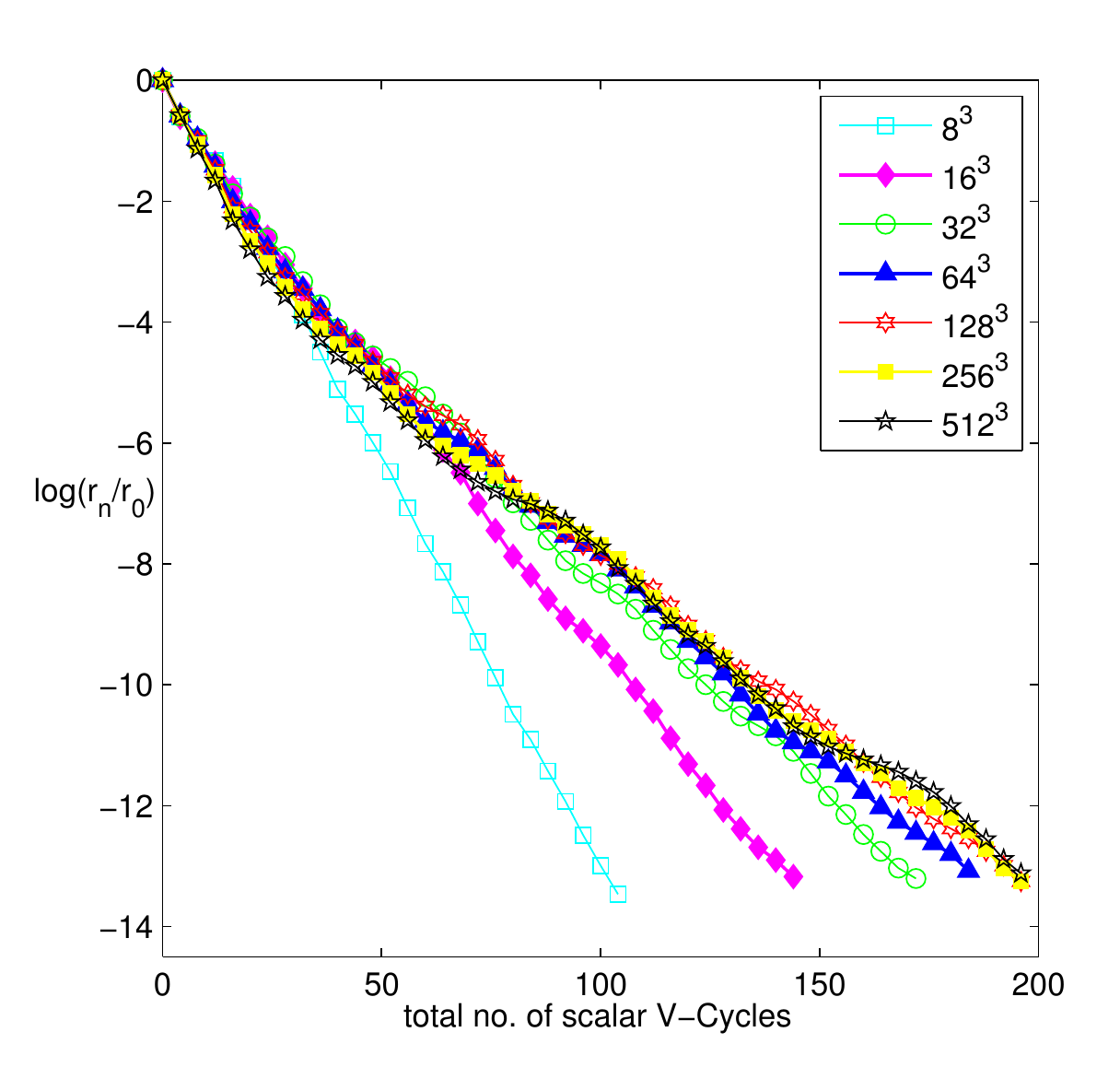}\includegraphics[width=0.49\textwidth]{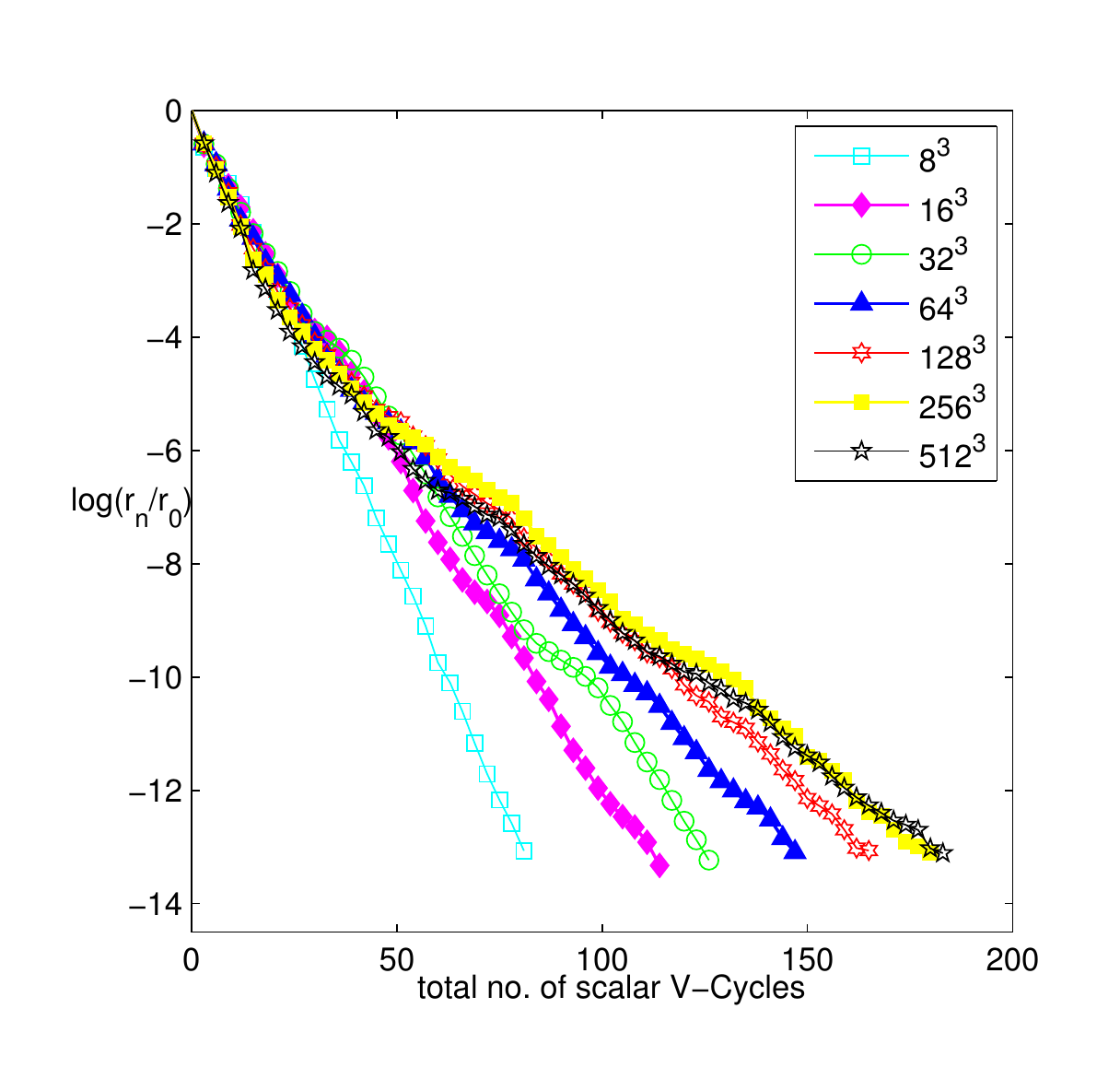}
\par\end{centering}

\caption{\label{fig:Scaling_100} The GMRES convergence history for preconditioners
$\V P_{1}$ (left panels) and $\V P_{2}$ (right panels) in two dimensions
(top panels) and in three dimensions (bottom panels) for the steady-state
bubble test problem with contrast ratio $r_{\mu}=r_{\rho}=100$, as
the grid size is varied.}
\end{figure}

An important goal in designing solvers suitable for large-scale calculations
is to ensure that the total number of multigrid cycles remains essentially
independent of the system size (or, equivalently, under grid refinement).
In Fig. \ref{fig:Scaling_100} we show convergence histories of GMRES
for varying grid sizes for the steady state bubble problem in both
two and three dimensions. In the left panels we show results for $\V P_{1}$
and the right panels for $\V P_{2}$. For this challenging variable-viscosity
problem (recall that the viscosity and density contrast ratio is $r_{\mu}=r_{\rho}=100$),
$\V P_{1}$ shows robust convergence for all of the grid sizes tested
here in both two and three dimensions, requiring no more than 200
multigrid V cycles (i.e., no more than $200/4=50$ GMRES iterations)
to reduce the residual to essentially roundoff tolerance even for
a $512^{3}$ grid. The convergence for preconditioner $\V P_{2}$
shows a very mild deterioration with increasing system size, although
the overall efficiency is still somewhat higher than $\V P_{1}$ for
all system sizes tested here. We have confirmed that making the subsolvers
nearly exact does not aid the overall GMRES convergence, despite the
substantial increase in the computational cost (data not shown).

\begin{figure}[!h]
\begin{centering}
\includegraphics[width=0.49\textwidth]{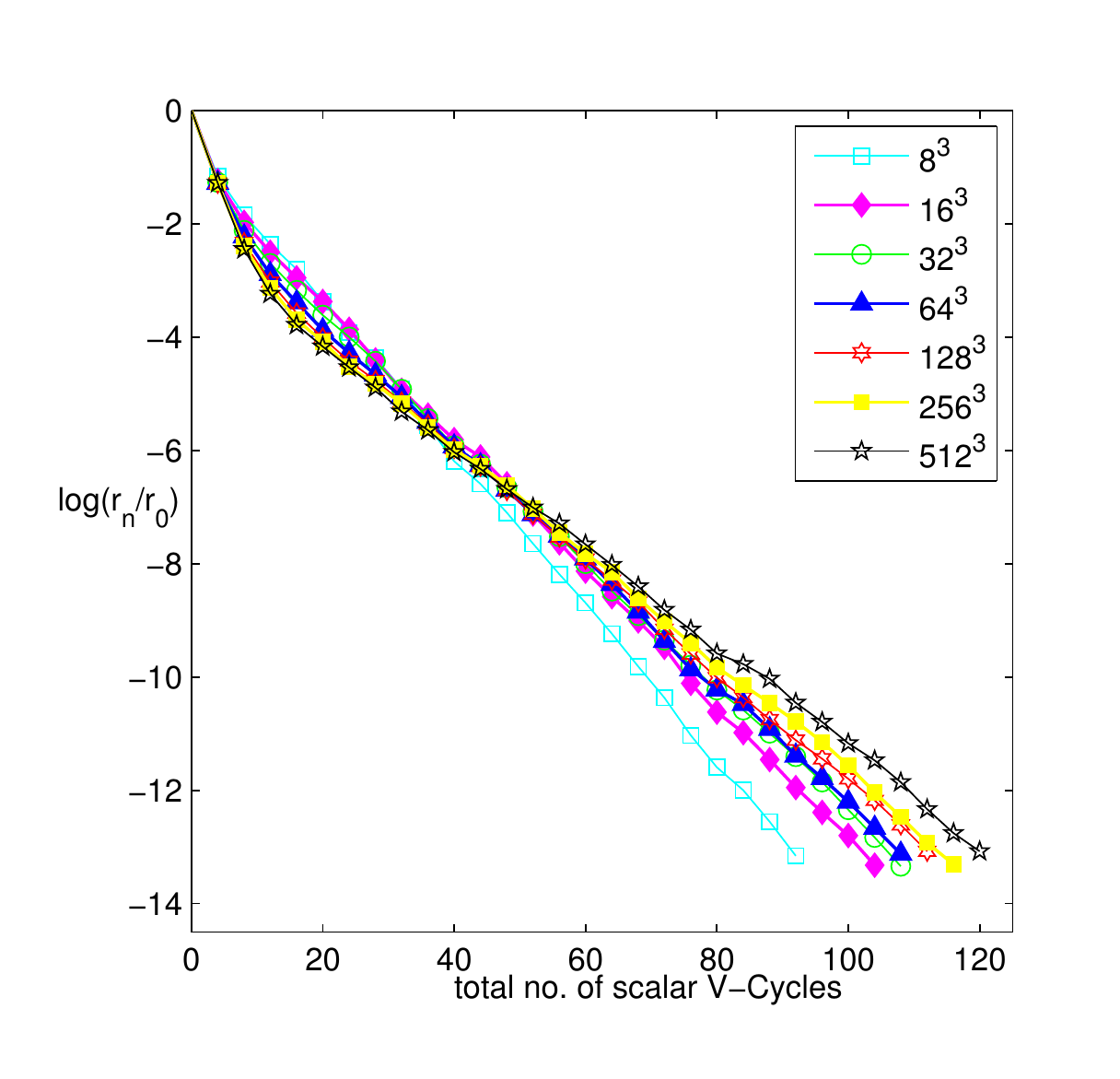}\includegraphics[width=0.49\textwidth]{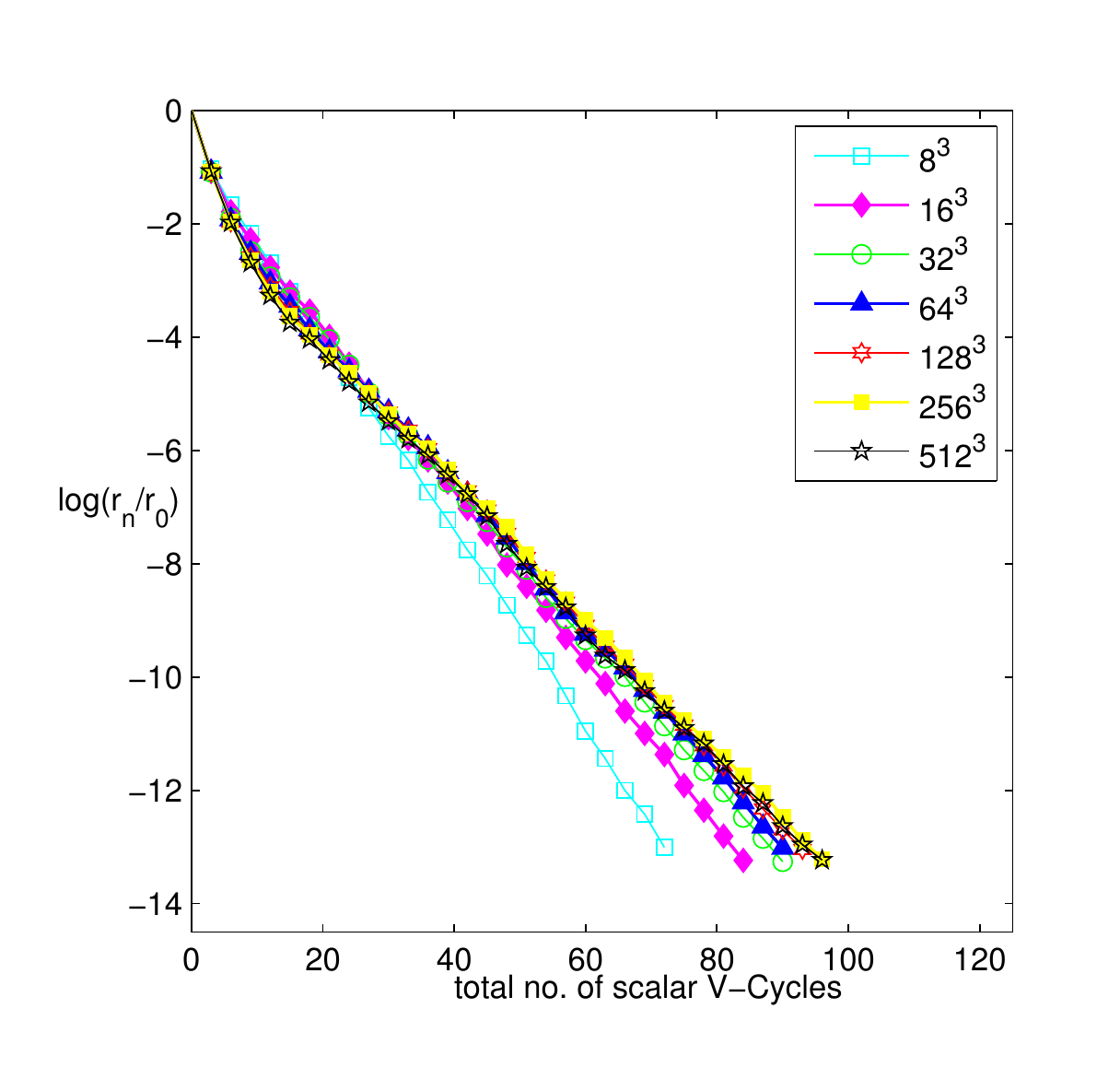}
\par\end{centering}

\caption{\label{fig:Scaling_2}Same as bottom two panels of Fig. \ref{fig:Scaling_100}
but for contrast ratio $r_{\mu}=r_{\rho}=2$.}
\end{figure}

It is important to point out that the exact convergence and its behavior
on system size depends sensitively on the details of the multigrid
algorithm (e.g., how the bottom level of the multigrid hierarchy is
handled, which is typically affected by parallelization), the restart
interval (here set to 10 iterations), and, most importantly, on the
contrast ratio. In Fig. \ref{fig:Scaling_2} we show scaling results
in three dimensions for a much weaker contrast ratio $r_{\mu}=r_{\rho}=2$.
In this case we see little to no effect of the system size on the
convergence rate, and the total number of GMRES iterations is less
than 30.

\section{\label{sec:Conclusions}Conclusions}

We studied several preconditioners for solving time-dependent and
steady discrete Stokes problems arising when solving fluid flow problems
on a staggered finite-volume grid. All of the preconditioners we studied
are based on approximating the inverse of the Schur complement with
a simple local operator and have been proposed before, though often
limited to either constant coefficient or steady flow. By suitably
approximating the inverse of the Schur complement in the case of time-dependent
variable-viscosity flow we were able to easily generalize these preconditioners
and thus substantially enlarge their practical applicability. Herein,
we modified and extended a previously proposed projection-based preconditioner
$\M P_{1}$ to variable-coefficient flows \cite{NonProjection_Griffith},
we generalized a well-known lower triangular preconditioner $\M P_{2}$
to variable-coefficient flow, and we extended a previously-studied
``fully coupled'' solver with a ``local viscosity'' preconditioner
\cite{LargeViscosityJump} to time-dependent flows to obtain an upper
triangular preconditioner $\M P_{3}$. The preconditioners investigated
here can be generalized to other stable or stabilized spatial discretizations
of the time-dependent Stokes equations, such as finite-element schemes
or adaptive mesh finite-volume discretizations.

Our primary focus was on studying the performance of these preconditioners
when the pressure and velocity subsolvers are performed on a uniform
staggered grid using geometric multigrid algorithms. We showed that
optimal convergence rates of the GMRES Stokes solver were obtained
when a single multigrid V cycle is employed as an inexact subdomain
solver. We numerically observed that all three preconditioners are
effective for both time-dependent and steady flow problems, with the
lower and upper triangular preconditioners being more efficient for
steady problems and $\M P_{1}$ being somewhat more efficient for
time-dependent problems, which is consistent with the findings of
Ref. \cite{NonProjection_Griffith}. All three preconditioners were
found to handle variable-coefficient problems rather well, with little
deterioration in convergence from the case of constant-coefficient
problems. Our observations are consistent with the conclusion of the
authors of Ref. \cite{LargeViscosityJump}, who ``find that it is
advantageous to use the FC {[}fully-coupled{]} approach utilizing
relaxed tolerances for solution of the sub-problems, combined with
the LV {[}local viscosity{]} preconditioner.'' 

All of our empirical observations are consistent with the general
observation that solving the coupled Stokes problem is comparable
to a single step of a fractional step method, and not more than 2-3
times more expensive than a fractional step even for difficult cases
of large contrast ratio, large viscous CFL number and non-trivial
boundary conditions. We believe that this mild increase in cost is
more than justified given the important advantages of the the coupled
approach when solving the Navier-Stokes equations. Furthermore, we
observed robust behavior of the projection and lower triangular preconditioners
for large systems with relatively frequent restarts. This demonstrates
that GMRES with preconditioners $\M P_{1}$ and $\M P_{2}$ provides
a robust solver for large-scale computations. In future studies, the
robustness of these preconditioners with respect to the variability
of viscosity and density should be studied more carefully. One aspect
of this is whether the spectrum of the preconditioned operator can
be provably bounded for arbitrary contrast ratios. More importantly,
however, the performance of the preconditioners in practical applications,
should be accessed. Experience with steady Stokes geodynamics applications,
which have extreme viscosity contrasts, are very promising \cite{LargeViscosityJump}.

\subsection*{Acknowledgments}

We thank Howard Elman for informative discussions. A. Donev and M.
Cai were supported in part by the NSF under grant DMS-1115341. Additional
support for A. Donev was provided by the DOE Office of Science through
Early Career award DE-SC0008271. B. E. Griffith acknowledges research
support from the National Science Foundation under awards OCI 1047734
and DMS 1016554. J. Bell and A. Nonaka were supported by the DOE Applied
Mathematics Program of the DOE Office of Advanced Scientific Computing
Research under the U.S. Department of Energy under contract No. DE-AC02-05CH11231.

\appendix

\section{\label{AppendixFourier}Fourier Analysis of Schur Complement}

The most important element in the preconditioners we study here is
the approximation of the Schur complement inverse. In previous work,
Fourier analysis \cite{CahouetChabard}, operator mapping properties
and PDE theory in \cite{MardalWinther2004,MardalWinther2011}, and
commutator properties (\ref{Schur_commutApp}) \cite{Kay,ApproximateCommutators,NonProjection_Griffith}
have been used to justify approximations to the Schur complement inverse.
Here we use Fourier analysis to justify our approximation to the Schur
complement inverse for the stress form of the viscous operator (\ref{div_stress}).
This analysis assumes periodic boundaries but should also inform the
case with physical boundary conditions.

For simplicity, we use two dimensional steady state Stokes equations
for illustration but extensions to three dimensions and time-dependent
flow are trivial. We denote the discrete Fourier transform of velocity
as $\hat{\V v}=\left[\hat{u},\,\hat{v}\right]^{T}$, and denote the
(purely imaginary) Fourier symbol of the staggered divergence operator
as $\hat{\Db}=[\hat{\Db}_{x}~\hat{\Db}_{y}]$, where $\Db_{x}$ and
$\Db_{y}$ represent the staggered finite difference operator along
the $x$ and $y$ axes. The Fourier transform of the staggered gradient
operator is $\hat{\Gb}=-\hat{\Db}^{\star}=[\hat{\Db}_{x}~\hat{\Db}_{y}]^{T}$,
and similarly, 
\[
\hat{L}_{p}=\hat{\Db}\hat{\Gb}=\hat{\Db}_{x}^{2}+\hat{\Db}_{y}^{2}=-\left(\abs{\hat{\Db}_{x}}^{2}+\abs{\hat{\Db}_{y}}^{2}\right).
\]

Our goal is to approximate the Schur complement inverse with a Laplacian-like
local operator $\V L_{\V S}$, i.e., to find $(\Db\Lb_{{\V{\mu}}}^{-1}\Gb)^{-1}=\V L_{\V S}.$
This is only an approximation in general but should be exact for periodic
constant-coefficient problems. In Fourier space, 
\begin{equation}
{\displaystyle \hat{\V L}_{\V S}=\left(\hat{\Db}\hat{\Lb}_{\V{\mu}}^{-1}\hat{\Gb}\right)^{-1}.}\label{Fourier_Schur}
\end{equation}
When the Laplacian form of the viscous term is used, $\Lb_{{\V{\mu}}}=\mu_{0}\Lb$,
we have 
\[
{\displaystyle \hat{\Lb}_{\V{\mu}}=\mu_{0}\left[\begin{array}{cc}
\hat{\Db}_{x}^{2}+\hat{\Db}_{y}^{2} & 0\\
0 & \hat{\Db}_{x}^{2}+\hat{\Db}_{y}^{2}
\end{array}\right]},
\]
which combined with (\ref{Fourier_Schur}) gives $\hat{\V L}_{\V S}=\mu_{0}.$
Applying an inverse Fourier transform gives the well-known result
$\V L_{\V S}=\mu_{0}\V I.$

When $\Lb_{{\V{\mu}}}$ is the discrete operator for the stress tensor
form of the viscous term (\ref{div_stress}) and the viscosity is
constant, we have 
\[
{\displaystyle \hat{\Lb}_{\V{\mu}}=\mu_{0}\left[\begin{array}{cc}
2\hat{\Db}_{x}^{2}+\hat{\Db}_{y}^{2} & \hat{\Db}_{y}\hat{\Db}_{x}\\
\hat{\Db}_{x}\hat{\Db}_{y} & \hat{\Db}_{x}^{2}+2\hat{\Db}_{y}^{2}
\end{array}\right]},
\]
which gives $\hat{\V L}_{\V S}=2\mu_{0},$ and therefore $\V L_{\V S}=2\mu_{0}\V I$.
This motivates our variable-viscosity generalization (\ref{Schur_app_varDen_varVis}).
When $\Lb_{{\V{\mu}}}$ is the discrete operator for the viscous term
with bulk viscosity and assuming both the shear viscosity and the
bulk viscosities are constant, we have 
\[
{\displaystyle \hat{\Lb}_{\V{\mu}}=\left[\begin{array}{cc}
(\frac{4}{3}\mu_{0}+\gamma_{0})\hat{\Db}_{x}^{2}+\mu_{0}\hat{\Db}_{y}^{2} & (\frac{1}{3}\mu_{0}+\gamma_{0})\hat{\Db}_{y}\hat{\Db}_{x}\\
(\frac{1}{3}\mu_{0}+\gamma_{0})\hat{\Db}_{x}\hat{\Db}_{y} & \gamma_{0}\hat{\Db}_{x}^{2}+(\frac{4}{3}\mu_{0}+\gamma_{0})\hat{\Db}_{y}^{2}
\end{array}\right],}
\]
which gives $\hat{\V L}_{\V S}=\left(\frac{4}{3}\mu_{0}+\gamma_{0}\right)$
and therefore ${\displaystyle \V L_{\V S}=\left(\frac{4}{3}\mu_{0}+\gamma_{0}\right)\V I}$.
This motivates our variable-viscosity generalization (\ref{eq:Schur_app_bulk}).

\section{\label{sec:AnalysisEigenvalues}Analysis of preconditioners with
exact subsolvers}

In this Appendix we give some analysis of the spectrum of the preconditioned
operators when exact pressure and velocity subsolvers are used. To
see how well the different preconditioners approximate the original
saddle point form (\ref{Saddle_sys}), we formally calculate 
\begin{eqnarray}
\V P_{1}^{-1}\V M & = & \left(\begin{array}{cc}
\Ib & (\Ib-\V{\rho}^{-1}\Gb\M L_{\V{\rho}}^{-1}\Db)\Ab^{-1}\Gb\\
\V 0 & \mathcal{S}^{-1}\V S
\end{array}\right),\label{preconditioned_P1M}
\end{eqnarray}
\begin{eqnarray}
\V P_{2}^{-1}\V M & = & \left(\begin{array}{cc}
\Ib & \Ab^{-1}\Gb\\
\V 0 & \mathcal{S}^{-1}\V S
\end{array}\right),\label{preconditioned_P2M}
\end{eqnarray}
\begin{eqnarray}
\V P_{3}^{-1}\V M & = & \left(\begin{array}{cc}
\Ib-\Ab^{-1}\Gb\mathcal{S}^{-1}\Db & \Ab^{-1}\Gb\\
\mathcal{S}^{-1}\Db & \V 0
\end{array}\right),\label{preconditioned_P3M}
\end{eqnarray}
and lastly
\begin{eqnarray}
\V P_{4}^{-1}\V M & = & \left(\begin{array}{cc}
\Ib & \Ab^{-1}\Gb\\
\mathcal{S}^{-1}\Db & \V 0
\end{array}\right).\label{preconditioned_P4M}
\end{eqnarray}
Recall that for constant-coefficient problems with exact subsolvers,
$\mathcal{S}^{-1}=-\theta\rho_{0}\M L_{p}^{-1}+\mu_{0}\V I$. For
periodic domains, the finite-difference operators $\Gb$, $\Db$,
$\M L$ and $\M L_{p}$ commute, 
\begin{equation}
\Gb\M L_{p}=\M L\Gb\quad\mbox{and}\quad\M L_{p}\Db=\Db\M L,\label{eq:commuting_periodic}
\end{equation}
and therefore $\V P_{1}^{-1}\V M$ is exactly the discrete identity
operator, and similarly, the (1,1) diagonal block of $\V P_{3}^{-1}\V M$
is zero.

From (\ref{preconditioned_P1M}), we see that the preconditioned system
is block upper triangular. Therefore, the eigenvalues of the preconditioned
system are either unity or the eigenvalues of $\mathcal{S}^{-1}\V S$.
Similarly, we can derive the eigenvalues of the preconditioned system
using (\ref{preconditioned_P2M}) and (\ref{preconditioned_P3M}).
Alternatively, one can write down the generalized eigenvalue system,
for instance, 
\[
\V M\left(\begin{array}{c}
\ub\\
p
\end{array}\right)=\lambda\V P_{3}\left(\begin{array}{c}
\ub\\
p
\end{array}\right)
\]
Again, one can see that the eigenvalues are either 1 or the eigenvalues
of $\mathcal{S}^{-1}\V S$.

When $\mu_{0}=0$, or equivalently, $\Delta t\rightarrow0$, we have
that $\V P_{1}^{-1}\V M=\Ib$, regardless of the boundary conditions.
If $\V P_{2}$ is used, we have 
\begin{eqnarray}
\V P_{2}^{-1}\V M & = & \left(\begin{array}{cc}
\Ib & \frac{\Delta t}{\rho_{0}}\Gb\\
\V 0 & \Ib
\end{array}\right),\label{preconditioned_ell_P2}
\end{eqnarray}
and therefore $(\V P_{2}^{-1}\V M-\Ib)^{2}=\V 0$. This proves that
in the inviscid case, the GMRES algorithm converges in 1 iteration
when preconditioner $\V P_{1}$ is used, and in 2 iterations when
$\V P_{2}$ or $\V P_{3}$ are used. When inexact subsolvers are used
our numerical results showed that all three preconditioners exhibit
exactly the same convergence rate in the inviscid case.

Furthermore, for constant viscosity ($\mu=1$) steady state ($\theta=0$)
problems on a two-dimensional domain of $n_{x}\times n_{y}$ grid
cells with no-slip boundaries, one can prove the following property
for the eigenvalues of the Schur complement $\V S=\Db\M L^{-1}\Gb$:
\begin{enumerate}
\item $\lambda(\V S)\in\{0\}\cup[\eta^{2},1]$, where $\eta$ is the inf-sup
constant independent of grid size \cite{Elman_book,Olshanskii}.
\item The multiplicity of the 0 eigenvalue is 1.
\item There are at most $2(n_{x}-1)+2(n_{y}-1)$ non-unit eigenvalues of
$\V S$.
\end{enumerate}
This is a quantitative statement of the intuitive expectation that
a few cells away from the boundaries $\V S$ is close to an identity
operator, just as for a periodic system (see Eq. (\ref{eq:commuting_periodic})).
The lower bound of the eigenvalues is a consequence of the uniform
div-stability \cite{Elman,Elman_book,Verfurth}. From (\ref{preconditioned_P1M})
(and also (\ref{preconditioned_P2M}) and (\ref{preconditioned_P3M})),
we see that the same conclusions hold for the preconditioned systems.
This analysis explains the good performance of the simple Schur complement
approximation even in the presence of nontrivial boundary conditions
\cite{NonProjection_Griffith}. Below we numerically compute the spectrum
of the eigenvalues of $\mathcal{S}^{-1}\V S$ for both constant and
variable viscosity steady flow, and confirm the theoretical predictions
above.

\subsection{Spectrum of the Preconditioned Operator}

\begin{figure}[!t]
\begin{centering}
\includegraphics[width=0.49\textwidth]{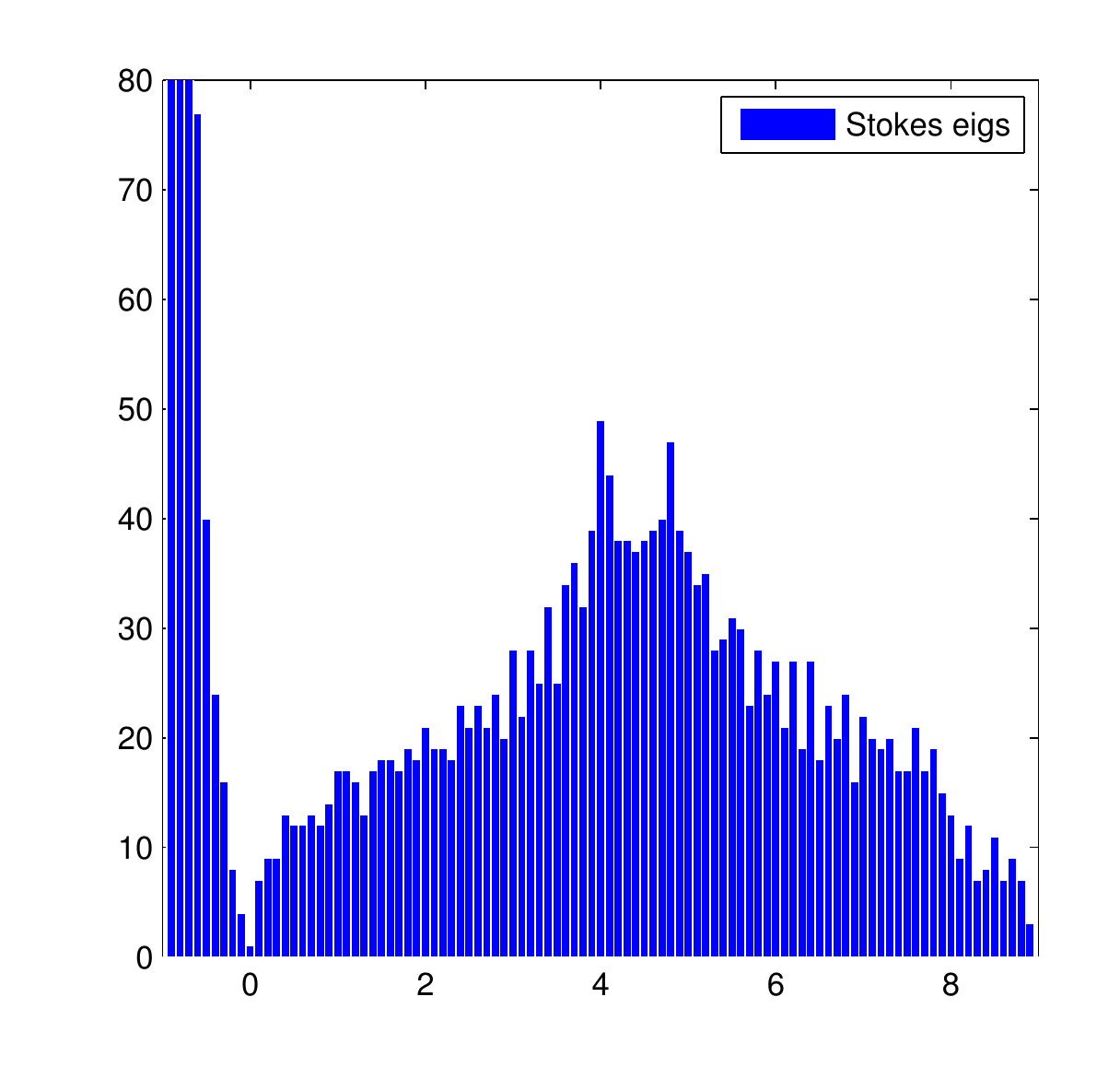}\includegraphics[width=0.49\textwidth]{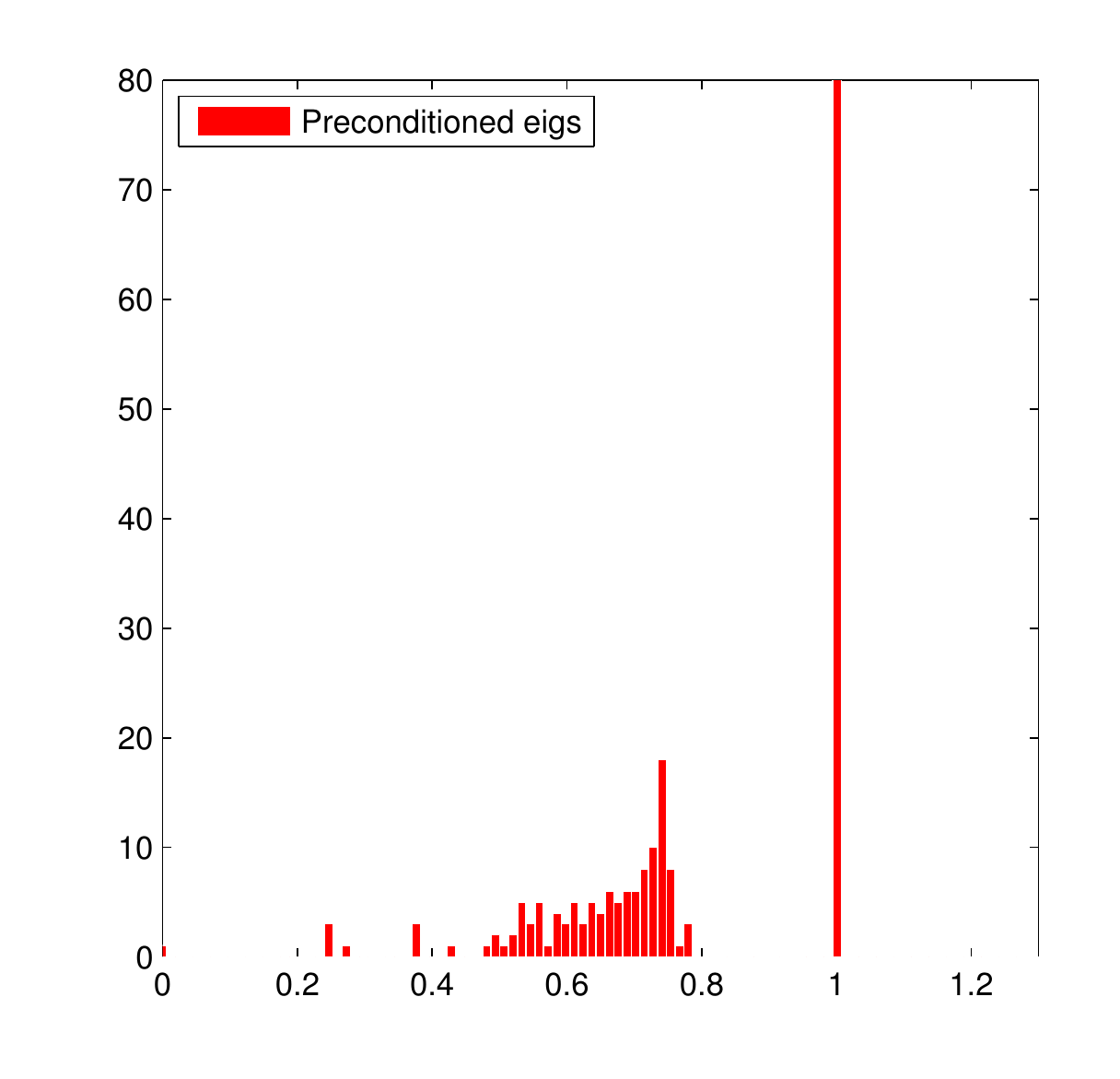}
\par\end{centering}

\begin{centering}
\includegraphics[width=0.49\textwidth]{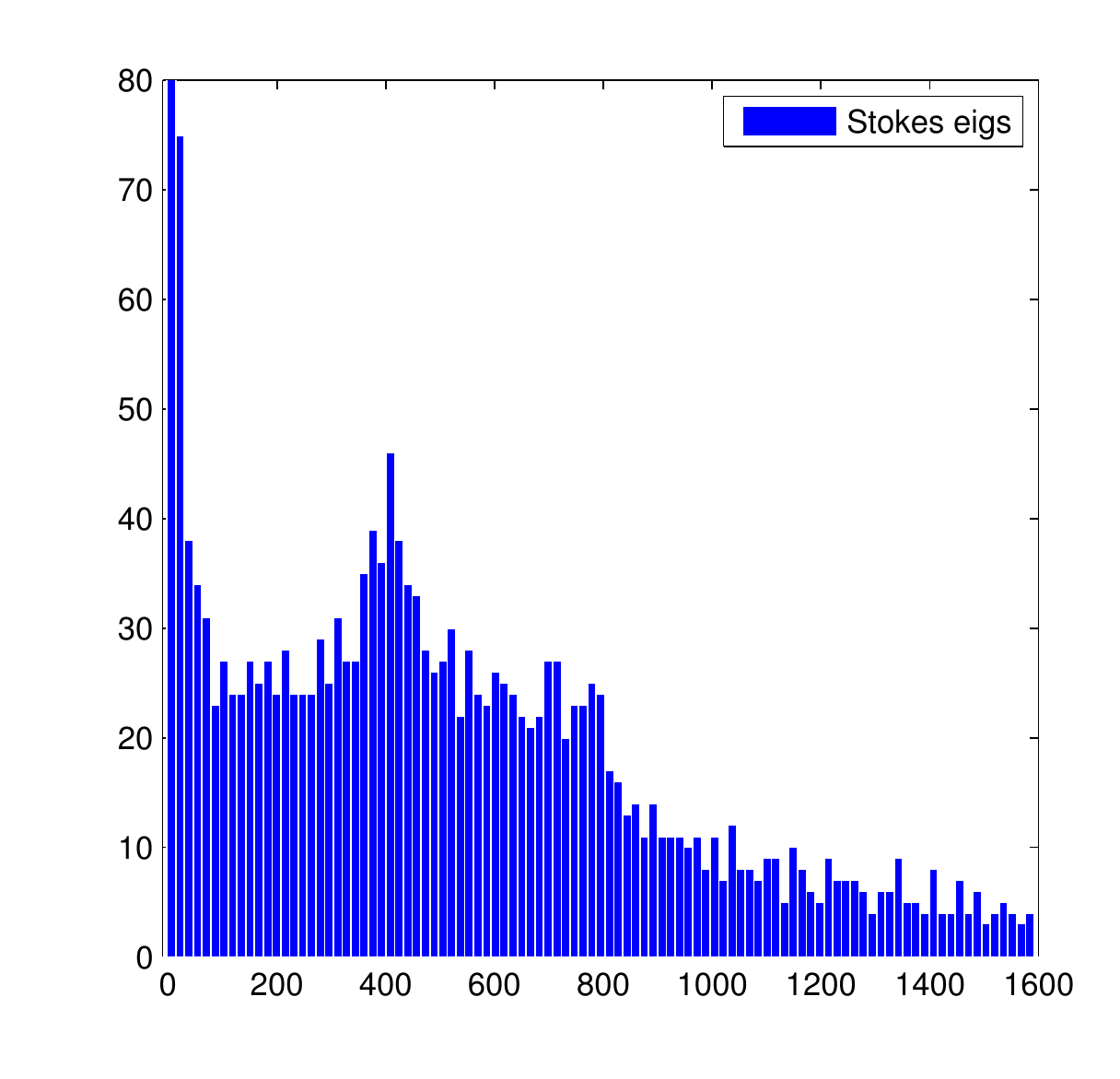}\includegraphics[width=0.49\textwidth]{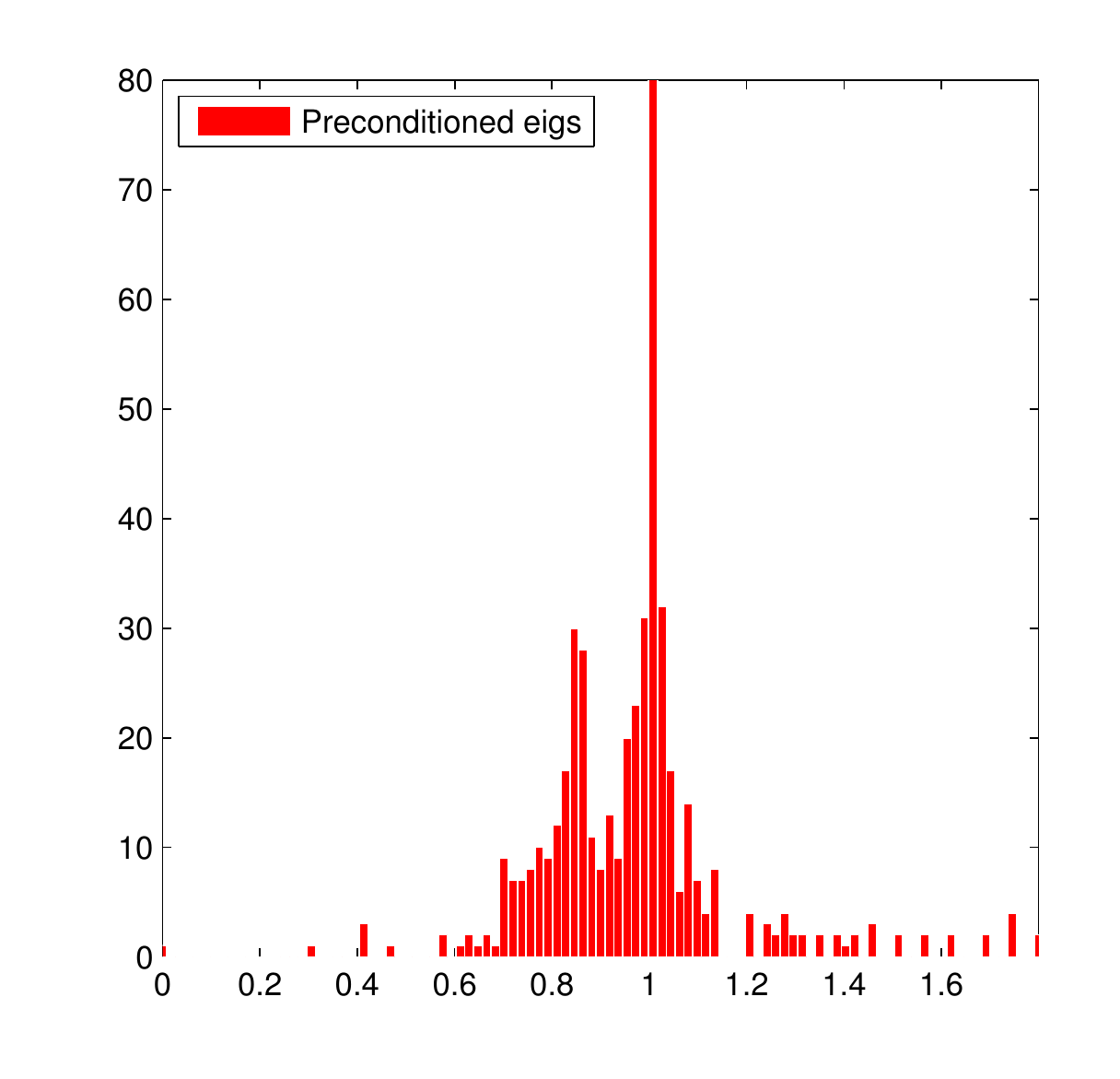}
\par\end{centering}

\caption{\label{fig:Eigenvalues} Histogram of the eigenvalues of $\M M$ (left
panels) and $\M P^{-1}\M M$ (more precisely, $\mathcal{S}^{-1}\V S$)
(right panels) for a steady Stokes problem on a grid of $32\times32$
cells with four no-slip boundaries. The vertical axis gives the number
of eigenvalues, truncated for the eigenvalue $\lambda=1$ in the right
panels due to the large number of unit eigenvalues. The case of constant
viscosity is shown in the upper panels, and the case of variable viscosity
with contrast ratio $100$ is shown in the lower panels.}
\end{figure}

Convergence analysis of the preconditioned GMRES method is not straightforward
and there is no simple link to the spectrum of the eigenvalues. Nevertheless,
it is generally accepted and widely observed that having closely clustered
eigenvalues of the preconditioned operator $\M P^{-1}\M M$ leads
to faster convergence. Furthermore, the ratio of the largest to the
smallest eigenvalue (excluding the trivial zero eigenvalue arising
from the fact pressure is only determined up to a constant) should
be bounded from above by a constant essentially independent of grid
size, and, possibly, viscosity and density contrast ratio.

We focus on the steady-flow case $\theta=0$ in two dimensions, for
a square domain of $n_{c}^{2}$ cells with four no-slip boundaries.
We consider using exact subdomain solvers, $\widetilde{\Ab}^{-1}=\M A^{-1}$
and $\widetilde{\V L}_{\V{\rho}}^{-1}=\M L_{\rho}^{-1}$ instead of
multigrid, relying on the fact that a well-designed multigrid cycle
is (essentially) spectrally-equivalent to an exact solver. We explicitly
form $\M M$ and $\mathcal{S}^{-1}\V S$ in MATLAB as dense matrices
and compute their eigenvalues. 

In Fig. \ref{fig:Eigenvalues}, we show a histogram of the eigenvalues
of the unpreconditioned and preconditioned operators for a square
domain of length $n_{c}=32$ cells with no-slip boundaries. In the
upper row in Fig. \ref{fig:Eigenvalues} we study the constant viscosity
case. The total number of DOFs is $N_{\text{dof}}=n_{c}^{2}+2n_{c}\left(n_{c}-1\right)=3008$.
Since the original Stokes system is of saddle point type, $\M M$
has both positive eigenvalues and negative eigenvalues, and there
are $n_{c}^{2}=1024$ eigenvalues that are smaller than or equal to
zero. While the unpreconditioned spectrum shows a broad spectrum of
eigenvalues with conditioning number that grows with the grid size,
the preconditioned spectrum shows that most eigenvalues are unity,
with the remaining $4\left(n_{c}-1\right)/N_{\text{dof}}\approx4\%$
nonzero eigenvalues remaining well-clustered.

In the lower row in Fig. \ref{fig:Eigenvalues} we study the variable
viscosity case for the bubble problem with viscosity contrast ratio
$r_{\mu}=100$. The unpreconditioned system is seen to be very badly
conditioned, with a broad spectrum of eigenvalues. By contrast, the
preconditioned operator is well-conditioned, with around 87\% of the
eigenvalues in the interval $\left(0.99,1.01\right)$. While some
eigenvalues are larger than unity in this case, the spread in the
eigenvalues is not much different from the constant-coefficient case.
This suggests that the spectrum remains localized around unity and
bounded away from zero even for rather large contrast ratios. It may
be possible to extend the finite-element theory developed in Refs.
\cite{OlshanskiiGrinevich,Grinevich} to prove that $\mathcal{S}^{-1}\V S$
is spectrally-equivalent to the identity matrix for the staggered
grid discretization we employ here.

\section{\label{sec:Multigrid}Multigrid algorithms}

We employ a standard V-cycle multigrid approach \cite{Briggs} for
both the cell-centered multigrid subsolver for the weighted Poisson
operator $\V L_{\V{\rho}}$ and the staggered velocity multigrid subsolver
for the viscous operator $\M L_{\V{\mu}}$. We use the standard residual
formulation, so that on all coarsened levels we are solving for the
error in the coarsened residual from the next-finer level. In our
implementation, the multigrid coarsening factor is 2, and coarsening
continues until the problem domain represented on the coarsest grid
contains two grid points (with respect to cell-centers) in any given
spatial direction. At the coarsest level of the multigrid hierarchy,
we perform a large (8 or more) number of relaxations, to ensure that
the preconditioner is a constant linear operator.

\begin{figure}
\centering{}\includegraphics[width=0.45\textwidth]{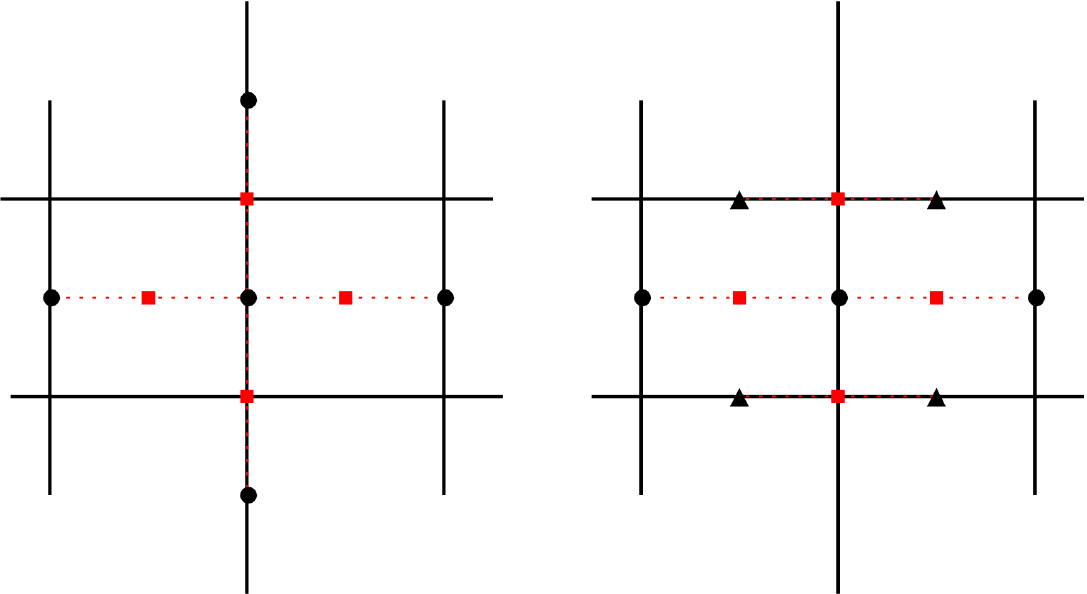}\hspace{1cm}\includegraphics[width=0.45\textwidth]{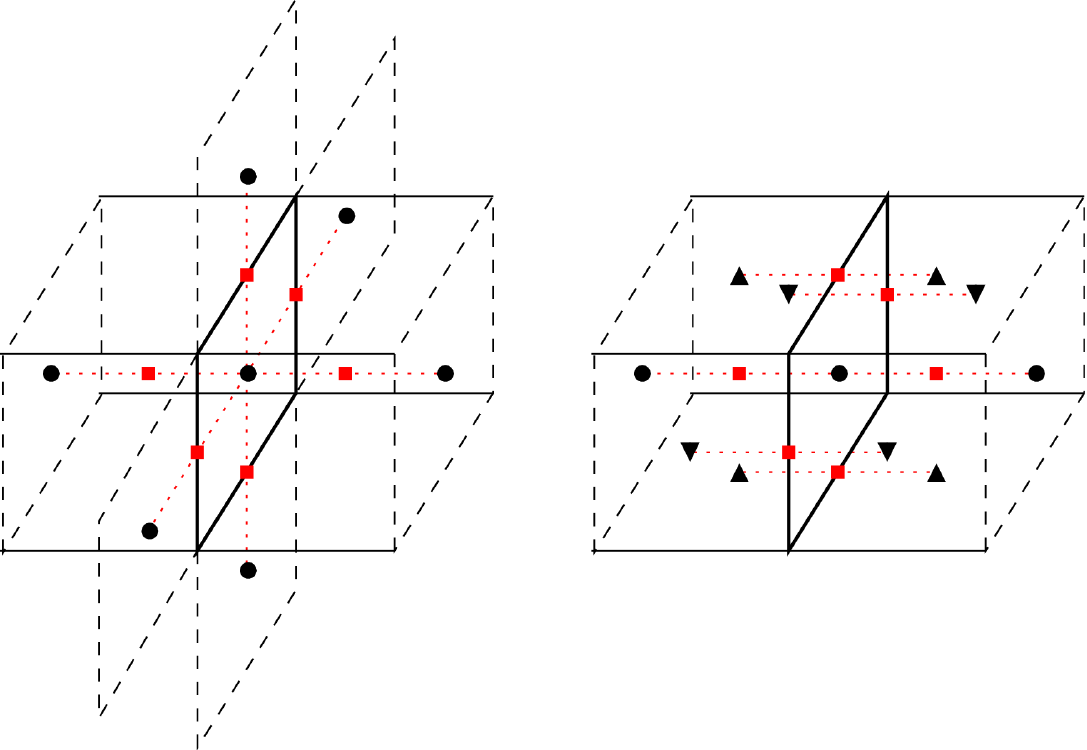}
\caption{\label{fig:viscOp}(\emph{Left panel}) The stencils for the $x$-component
of $\grad\cdot\beta\grad\phib$ (left) and $\grad\cdot\beta(\grad\phib)^{T}$
(right), in two dimensions. The black circles indicate locations of
$u$. The black triangles indicate locations of $v$. The red dots
indicate the location of the $\beta$ and the gradients of velocity.
(\emph{Right panel}) The stencils for the $x$-component of $\grad\cdot\beta\grad\V v$
(left) and $\grad\cdot\beta(\grad\V v)^{T}$ (right), in three dimensions.
The black circles indicate locations of $u$. The black triangles
indicate locations of $v$ and $w$. The red dots indicate the location
of the $\beta$ and the gradients of velocity.}
\end{figure}

Multigrid consists of 3 major components: (i) choice of relaxation
at a particular level, (ii) coarsening/restriction operator, and (iii)
interpolation/prolongation operator.

\textbf{Relaxation.} Both the staggered and cell-centered solvers
use multicolored Gauss-Seidel smoothing. The cell-centered solver
uses standard red-black relaxation, whereas the staggered solver uses
a $2d$-colored relaxation, where $d$ is the dimensionality of the
problem. Because the coupling between the degrees of freedom corresponding
to a given component of velocity is the same as for the cell-centered
Poisson equation, by coloring each component of the velocity separately,
as in the standard red-black coloring (i.e., coloring odd grid points
with a different color from the even grid points), we obtain decoupling
between the $2d$ colors so that each color can be relaxed separately
(improving convergence and aiding parallelization). We relax the components
of velocity in turn (i.e., in three dimensions, we order the relaxations
as red-$x$, black-$x$, red-$y$, black-$y$, red-$z$, black-$z$),
although other orderings of the colors are possible. Refer to Figure
\ref{fig:viscOp} for a physical representation of the viscous operator
stencils. Given a cell-centered operator of the form, $\grad\cdot\beta\grad\phi\equiv\mathcal{L}\phi=r$,
or a staggered operator of the form, $\alpha\phi-\grad\cdot\beta[\grad\phi+(\grad\phi)^{T}]\equiv\mathcal{L}\phi=r$,
the relaxation takes the form 
\begin{equation}
\phi^{k+1}=\phi^{k}+\omega\mathcal{D}^{-1}(r-\mathcal{L}\phi^{k}),
\end{equation}
for each color in turn, where the superscript represents the iterate,
and $\mathcal{D}^{-1}$ is the inverse of the diagonal elements of
$\mathcal{L}$. We use unit weighting factor%
\footnote{Note that for Jacobi relaxation with the stress form of the viscous
operator, a standard analysis suggests $\omega=1/2$ as the optimal
relaxation parameter (ensuring damping of all modes).%
}, $\omega=1$ (suggested to be near-optimal in numerical experiments)
for both subsolvers.

\textbf{Restriction.} For the cell-centered solver, restriction is
a simple averaging of the $2^{d}$ fine cells. For the staggered solver,
we use a slightly more complicated 6-point ($d=2$) or 12-point ($d=3$)
stencil. For example, for $x$-faces we use 
\begin{equation}
\phi_{i,j}^{{\rm c}}=\frac{1}{8}\left(\phi_{2i-1,2j}^{{\rm f}}+\phi_{2i-1,2j+1}^{{\rm f}}+\phi_{2i+1,2j}^{{\rm f}}+\phi_{2i+1,2j+1}^{{\rm f}}\right)+\frac{1}{4}\left(\phi_{2i,2j}^{{\rm f}}+\phi_{2i,2j+1}^{{\rm f}}\right)
\end{equation}
As seen in Figure \ref{fig:viscOp}, for the staggered solver we require
$\alpha$ at faces, and $\beta$ at both cell-centers and nodes ($d=2$)
or edges ($d=3$). When creating coefficients at coarser levels, we
obtain $\alpha$ by averaging the overlaying fine faces, cell-centered
$\beta$ by averaging the overlaying fine cell-centered values, $\beta$
on nodes through direct injection, and $\beta$ on edges by averaging
the overlaying fine edges.

\textbf{Prolongation.} For the cell-centered solver, prolongation
is simply direct injection from the coarse cell to the overlaying
$2^{d}$ fine cells. For the staggered solver, we use a slightly more
complicated stencil that involves linear interpolation for fine faces
that overlay coarse faces, and bilinear interpolation for fine faces
that do not overlay coarse faces. For example, for $x$-faces we use
\begin{equation}
\phi_{i,j}^{{\rm f}}=\frac{3}{4}\phi_{i/2,j/2}^{{\rm c}}+\frac{1}{4}\phi_{i/2,j/2-1}^{{\rm c}},\quad\mbox{for }i\mbox{ and }j\mbox{ both even},
\end{equation}
\begin{equation}
\phi_{i,j}^{{\rm f}}=\frac{3}{8}\left(\phi_{i/2,j/2}^{{\rm c}}+\phi_{i/2+1,j/2}^{{\rm c}}\right)+\frac{1}{8}\left(\phi_{i/2,j/2-1}^{{\rm c}}+\phi_{i/2+1,j/2-1}^{{\rm c}}\right),\quad\mbox{for }i\mbox{ odd and }j\mbox{ even},
\end{equation}
where we use integer division in the index subscripts.


\end{document}